\documentclass[12pt,
twoside,a4paper]{amsart}
\usepackage[T1]{fontenc}
\usepackage[latin1]{inputenc}
\usepackage[left=25mm,top=30mm,bottom=30mm]{geometry}
\usepackage{amssymb}
\usepackage{mathtools}
\usepackage{lmodern}
 \usepackage[
 final,
  scrtime
 ]{prelim2e}

  \renewcommand{\PrelimText}{\footnotesize[\,Version: 
  \texttt{\jobname.tex}\hfill \today\ at \thistime\,]}
      \theoremstyle{plain}
      \newtheorem{Thm}{Theorem}[section]
      \newtheorem{Lem}[Thm]{Lemma}
      
      \newtheorem{Rem}[Thm]{Remark}
      \theoremstyle{definition}
            
      \theoremstyle{remark}
\let\oldsqrt\sqrt
\def\sqrt{\mathpalette\DHLhksqrt}
\def\DHLhksqrt#1#2{%
\setbox0=\hbox{$#1\oldsqrt{#2\,}$}\dimen0=\ht0
\advance\dimen0-0.2\ht0
\setbox2=\hbox{\vrule height\ht0 depth -\dimen0}%
{\box0\lower0.4pt\box2}}
\newcommand{\N}{{\mathbb N}}     % 
\newcommand{\R}{{\mathbb R}}     %
\newcommand{\Z}{{\mathbb Z}}%
\renewcommand{\epsilon}{\varepsilon}\renewcommand{\phi}{\varphi} %
\renewcommand{\rho}{\varrho}\renewcommand{\theta}{\vartheta}     %
\DeclareMathOperator*{\bigconv}{\mbox{\LARGE$\ast$}}
\renewcommand{\[}{\begin{eqnarray*}}\renewcommand{\]}{\end{eqnarray*}}
\newcommand{\la}{\begin{eqnarray}}\newcommand{\al}{\end{eqnarray}}

\newcommand{\dd}{{\,\mathrm{d}}}
\begin{document}

% First we specify the top matter (author, title, etc).
%
% Note: All of the top matter items are optional and can be omitted.
% But you will probably want to specify at least the author and title
% and perhaps an abstract.

% author information

% first author 

\author[Mattner]{Lutz Mattner}
\address{Universit\"at Trier, Fachbereich IV -- Mathematik, 54286~Trier, Germany}
\email{\{mattner, schu4502\}@uni-trier.de}

% second author

\author[Schulz]{Jona Schulz}

%    % the address where the research was carried out
%  \address{Universit\"at Trier, Fachbereich IV -- Mathematik, 54286~Trier, Germany}

% current address, usually not needed because it is the same as the
% regular address
%    \curraddr{Department of Mathematics, Pennsylvania State University,
%      University Park, State College PA 16802}

%\email{schu...@uni-trier.de}
   
% title
\title[On normal approximations to symmetric hypergeometric laws]
{On normal approximations to\\ symmetric hypergeometric laws}

% Note that the short title for running heads goes in square
% brackets.  This is optional.  The long title goes in curly
% braces.  In the long title, line breaks are indicated by \\.

% abstract (optional)
\begin{abstract}
The Kolmogorov distances between a symmetric hypergeometric law 
with standard deviation $\sigma$ and its usual normal approximations
are computed and shown to be less than $1/(\sqrt{8\pi}\,\sigma)$, 
with the order $1/\sigma$ and the constant $1/\sqrt{8\pi}$ being optimal.  
The results of Hipp and Mattner~(2007) for symmetric binomial laws 
are obtained as special cases.

Connections to Berry-Esseen type results in more general situations concerning 
sums of simple random samples or Bernoulli convolutions are explained.

Auxiliary results of independent interest include rather sharp normal distribution 
function inequalities, a simple identifiability result for hypergeometric laws,
and some remarks related to L\'evy's concentration-variance inequality.   
\end{abstract}

% AMS subject classifications (used in AMS journals)

% AMS keywords (used in AMS journals)

% acknowledge support, etc
% \thanks{This research was partially supported by DFG grant MA 1386/3-1.}
% \thanks{This paper constitutes a part of the second author's diploma thesis.}

% dedication
%    \dedicatory{Dedicated to Professor Donald Knuth on the occasion
%      of his $100$th birthday}

% today's date, or fill in whatever date you prefer
 \date{April 24, 2014}

% This ends the top matter information.
% We can now tell LaTeX to display the top matter.

\maketitle

\tableofcontents

% Having displayed the top matter, we now proceed to the body of the
% article. The body of the article is divided into sections.
% Each section begins with a \section command.

\section{Introduction and main result}
\subsection{Aim}
This paper generalizes the error bound in the central limit theorem
for symmetric binomial laws of Hipp and Mattner~\cite{HM},
which up to now was the only nontrivial example of a 
Berry-Esseen type inequality with an optimal constant
known to the present authors,
to a still optimal bound covering also symmetric hypergeometric laws. 
These solutions of special cases of the Berry-Esseen problem
are of some particular interest for more general situations, as we attempt 
to explain in the subsection~\ref{Subsec:Background}   below,
and are also remarkable in view of the apparent 
difficulty of determining merely close to optimal
Berry-Esseen type inequalities in related special situations, as witnessed  
by the recent investigations of arbitrary binomial laws by 
Nagaev and Chebotarev~\cite{NagCheb} and of arbitrary Bernoulli 
convolutions, which include in particular all hypergeometric laws as is known from
\cite{Vat_Mik},  by Neammanee~\cite{Neammanee}.

\subsection{Background: Berry-Esseen  % theorems 
for sampling  with or without replacement}          \label{Subsec:Background}
Throughout this paper, let 
$\Phi$ denote the distribution function of the standard normal law.
In this subsection, let $g:\mathopen[1,\infty\mathclose[ \rightarrow 
\mathopen]0,\infty\mathclose]$ denote the pointwise smallest function such that 
\begin{eqnarray} 
 \left\| F -\Phi\right\|_\infty 
 &\le& \frac{g\left(\frac{\beta}{\sigma^3} \right)}{\sqrt{n}}
\end{eqnarray}
holds whenever $n\in\mathbb N$  and $F$ is the distribution function of the 
standardized sum of $n$ i.i.d.~random variables with law $P$
on the real line $\mathbb R$ with mean $\mu$, variance $\sigma^2>0$, and finite third
centred absolute moment $\beta = \int|x-\mu|^3\,\mathrm{d} P(x)$.
Let  further $C \in \mathopen]0,\infty\mathclose]$ denote the smallest
constant such that $g(\rho)\le C\rho$ holds for every $\rho\in \mathopen[1,\infty\mathclose[$.
Then the classical Berry-Esseen theorem for sums of i.i.d.~random variables states 
that $C<\infty$. More recent investigations aim, among other goals,
at obtaining rather sharp upper bounds on the function $g$, and here the best result
announced  so far appears to be Shevtsova's~\cite{Shev2013} bound
$ \min\{0.4690\rho, 0.3322(\rho + 0.429), 0.3031(\rho + 0.646)\}$
for each $\varrho$, which, when combined with a classical lower bound for $C$
due to Esseen~\cite{Ess1956}, yields in particular  
$0.4097<(\sqrt{10}+3)/(6\sqrt{2\pi})\leq C <0.4690$, 
and $g(1)<0.4690$. However, as,
by a discussion of equality in Lyapunov's moment inequality,
 $\beta/\sigma^3=1$ iff $P$ is a uniform law on 
two points, without loss of generality $0$ and $1$,
the special Berry-Esseen theorem for 
symmetric binomial laws~\cite[Corollary 1.2]{HM} yields 
$g(1)=1/\sqrt{2\pi}< 0.3990$. Although, unfortunately, we do not yet know whether
$g$ is continuous at $1$, the cited special result suggests the possibility 
of an improvement of  Shevtsova's bound for $\rho$ close to $1$.  

Analogously, the Berry-Esseen type theorem for sampling without replacement 
from a finite population due to H\"oglund~\cite{H} can be stated as follows:
Let  $h:\mathopen[1,\infty\mathclose[ \rightarrow 
\mathopen]0,\infty\mathclose]$ denote the pointwise smallest function such that
\begin{eqnarray}
  \left\| F -\Phi\right\|_\infty  
 &\le&  \frac{h\left(\frac{\beta}{\sigma^3}\right)}{\sqrt{n\,(1-\frac{n}{N})} }
\end{eqnarray}
holds whenever $N\in\mathbb{N}$ and $x\in\mathbb{R}^N$ are such that the law 
$P\coloneqq \frac1N\sum_{i=1}^N\delta_{x_i}$ has mean $\mu=\frac1N\sum_{i=1}^Nx_i$,
variance $\sigma^2= \frac1N\sum_{i=1}^N (x_i-\mu)^2 >0$, and the third
centred absolute moment $\beta=\frac1N\sum_{i=1}^N |x_i-\mu|^3 $, and 
whenever $n\in\{1,\ldots,N-1\}$ and $F$ is the distribution function of 
\begin{eqnarray}
    \frac{S-n\mu}{\sqrt{n\,(1-\frac{n}{N})}\sigma}
\end{eqnarray}
with $S$ being the sum of a simple random sample of size $n$ from $x$. 
Let  further $D \in \mathopen]0,\infty\mathclose]$ denote the smallest
constant such that $h(\rho)\le D\rho$ holds for every $\rho\in \mathopen[1,\infty\mathclose[$.
Then H\"oglund's theorem states that $D<\infty$. With $g$ and $C$ as in the previous paragraph, 
we have the simple Lemma~\ref{Lem:1.1} below, and hence $C \le D$, 
but  we are not aware of any published explicit upper bounds for $h$ or $D$.
However, using again that $\beta/\sigma^3 =1$ iff $P$ is a uniform law on 
two points,  we see that the special Berry-Esseen theorem
for symmetric hypergeometric laws~\ref{Thm:Main} below  and the formula for 
$\sigma_0$ in~\eqref{Eq:sigma,sigma_0}
(where $F$ and $\sigma$ have different meanings)  yield
\[
 h(1) &=& \sup\left\{ \sqrt{n\left(1-\tfrac{n}{N}\right)}\,d \,:\, 
            d,n,N \text{ as in Theorem~\ref{Thm:Main}(a)}  \right\} \\
   &=&\sup\left\{ 2\sigma_0 d \,:\, 
            d,\sigma_0 \text{ as in Theorem~\ref{Thm:Main}(a)}  \right\} 
   \,\ =\,\ \frac{1}{\sqrt{2\pi}} 
\]
by Remark~\ref{Rem:Complements_to_main_result}(b) with $\tau=\sigma_0$, and by using
the optimality of $\frac1{\sqrt{8\pi}}$ from  Theorem~\ref{Thm:Main}(a), or 
$g(1)=  \frac1{\sqrt{2\pi}}$  and Lemma~\ref{Lem:1.1}.
Hence $h(1)=g(1)$, suggesting that any effective upper bounds for $h(\rho)$ which might 
become available in the future should be close to $1/\sqrt{2\pi}$ for $\rho$ close to one, 
and perhaps even close to $g(\rho)$  in any case. 
Again, unfortunately, we do not yet know whether $h$ is continuous at $1$.

%%%%%%%%%%%%%%%%%%%%%%%%%%%%%%%%%%%%%%%%%%%%%%%%%%%%%%%%%%%%%  Lemma 1.1
\begin{Lem}                                                \label{Lem:1.1}
The functions $g$ and $h$ introduced above satisfy $g\le h$.
\end{Lem}
\begin{proof} Given $\rho\in[1,\infty[$ and any $\gamma\in\R$ with  
$\gamma < g(\rho)$, the definition of $g(\rho)$ as a supremum yields an $n\in\N$
and a law $P$ on $\R$ with third standardized absolute moment
$\beta/\sigma^3 =\rho$
and, using a reflection argument if necessary, an $s\in\R$
with $\Delta\coloneqq \sqrt{n}\left(P^{\ast n}\left(\mathopen]0,s\mathclose[\right) - 
\Phi\left( \frac{s-n\mu}{\sqrt{n}\sigma}\right)\right)>\gamma $.
Using the denseness with respect to weak convergence of the laws with finite support and rational 
point masses following from~\cite[Theorem 15.10]{AliprantisBorder}
together with a simple truncation argument, we can take  
$x_{N} \in\R^N $ for $N>n$ such that 
$P_N\coloneqq \frac1{N}\sum_{i=1}^N \delta_{x_{N,i}}$ converges to $P$ weakly 
and together with its moments and absolute moments up to the third order, for 
$N\rightarrow\infty$. 
Since the law $Q_N$ of the sum of a simple random sample of size 
$n$ from $P_N$ differs from $P_N^{\ast n}$ in the supremum distance
by at most $\frac{n(n-1)}{2N}$, see \cite{Freedman}, 
and since $P_N^{\ast n}$ tends weakly to 
$P^{\ast n}$ for $N\rightarrow \infty$, we get
\[
 \frac{h(\rho)}{ \sqrt{n} } 
   &=& \lim_{N\rightarrow\infty}\frac{h(\rho)}{\sqrt{ n\left(1- \tfrac{n}{N}\right)}}
      \,\ \ge \,\  \lim_{N\rightarrow\infty}\left(Q^{}_N\left(\mathopen]-\infty,s\mathclose[ \right)
     - \Phi\left( \frac{s-n\mu}{ \sqrt{ n\left(1- \frac{n}{N}\right)}\,\sigma}  \right) \right)
  \,\ \ge \,\ \frac{\Delta}{\sqrt{n}}
\]
and hence $h(\rho) > \gamma$.
\end{proof}

%%%%%%%%%%%%%%%%%%%%%%%%%%%%%%%%%%%%%%%%%%%%%%%%%%%%%%%%%%%%%%%%%%%%%%%%%%%%%%%%%%%%%%%%%%
%%%%%%%%%%%%%%%%%%%%%%%%%%%%%%%%%%%%%%%%%%%%%%%%%%%%%%%%%%%%%%%%%%%%%%%%%%%%%%%%%%%%%%%%%%
\subsection{Hypergeometric laws}
%%%%%%%%%%%%%%%%%%%%%%%%%%%%%%%%%%%%%%%%%%%%%%%%%%%%%%%%%%%%%%%%%%%%%%%%%%%%%%%%%%%%%%%%%%
%%%%%%%%%%%%%%%%%%%%%%%%%%%%%%%%%%%%%%%%%%%%%%%%%%%%%%%%%%%%%%%%%%%%%%%%%%%%%%%%%%%%%%%%%%
Let us here formally define hypergeometric and a few related laws on $\R$
and collect some standard properties of them. For $a\in\R$, we write $\delta_a$ for 
the {\em Dirac} measure concentrated at $a$.
For $\alpha\in\R$, we write $\alpha_{}^{\underline{k}} \coloneqq \prod_{j=1}^{k}(\alpha-j+1)$
and $\binom{\alpha}{k}\coloneqq \alpha_{}^{\underline{k}}/k!$ for $k\in\N_0$,
and, with the exception of the proof of Lemma~\ref{Lem:w(x)}, 
we put in this paper $\binom{\alpha}{k} \coloneqq  0$ if $k\notin\N_0$.
Then, for $n\in\N_0$ and $p\in[0,1]$, the binomial law $\mathrm{B}_{n,p}$ 
can be defined by $\mathrm{B}_{n,p}(\{k\}) \coloneqq \mathrm{b}_{n,p}(k)
\coloneqq  \binom{n}{k}p^k(1-p)^{n-k}$ 
for $k\in \Z$; and a law $P$ is {\em Bernoulli} if $P=\mathrm{B}_{1,p}$ for 
some $p\in[0,1]$.
For $r,b\in\N_0$ and $n\in\{0,\ldots,r+b\}$, we let $\mathrm{H}_{n,r,b}$ denote 
the {\em hypergeometric}
law of the number of red balls drawn in a simple random sample of size $n$ from 
an urn containing $r$ red and $b$ blue balls (red and blue, 
and not for example black and white, since the present choice 
of the colours leads to the same initial letters in several languages),
so that we have
\la                                              \label{Eq:H_nrb_densitiy}
 \mathrm{H}_{n,r,b}(\{k\}) &\eqqcolon  & \mathrm{h}_{n,r,b}(k)
  \,\ = \,\ \frac{\binom{r}{k}\binom{b}{n-k}}{\binom{r+b}{n}} 
     \qquad \text{ for }k\in\Z, 
\al
which may also be used to define $\mathrm{H}_{n,r,b}$   to avoid reference  
to a sampling model. No confusion of the notation $\mathrm{h}_{n,r,b}$ with the letter
$h$ used for various objects in this paper seems likely.
We use the convention $\frac00 \coloneqq 0$, relevant for example
in~\eqref{Eq:Hyp_mean_var_kappa} below if $r+b\in\{0,1,2\}$.
Except for the trivial cases of $n=0$ or $p=0$, a binomial law $\mathrm{B}_{n,p}$ uniquely determines
its parameters $n$ and $p$, and is symmetric about its mean iff $p=\frac12$, in which case 
the mean is $\frac{n}2$. The following lemma collects analogous or related simple facts for 
hypergeometric laws, used below but apparently not easily available 
from the literature.  

%%%%%%%%%%%%%%%%%%%%%%%%%%%%%%%%%%%%%%%%%%%%%%%%%%%%%%%%%%%%%%%%%%%%%%%%% Lemma 1.2
\begin{Lem}                                \label{Lem:Hyp_basics}
%%%%%%%%%%%%%%%%%%%%%%%%%%%%%%%%%%%%%%%%%%%%%%%%%%%%%%%%%%%%%%%%%%%%%%%%%%%%%%%%%%%%%%%%
Let $r,b\in\N_0$ and $n\in\{0,\ldots,r+b\}$. 

\smallskip\noindent{\rm\textbf{(a) Some basic descriptive properties.}}
$\mathrm{H}_{n,r,b}$ has the support
\la                                                  \label{Eq:Supp_Hyp} 
  \{k\in \Z: \mathrm{h}_{n,r,b}\left(k\right) >0\} 
    &=& \{ (n-b)_+ , \ldots , n\wedge r\}
\al
and the first three cumulants (mean, variance, third centred moment)
\la                                        \label{Eq:Hyp_mean_var_kappa}
 &&
 \mu = \frac{nr}{r+b} , 
  \quad \sigma^2=\frac{nrb(r+b-n)}{(r+b)^2(r+b-1)},
  \quad \kappa_3 = \frac{nrb(b-r)(r+b-n)(r+b-2n)}{(r+b)^3(r+b-1)(r+b-2)}.
\al

\smallskip\noindent{\rm\textbf{(b) (Non-)identifiability of parameters.}}
We have 
\la                             \label{Eq:Hyp_parameters_unidentifiable_1}  % (7)
 \mathrm{H}_{n,r,b}= \mathrm{H}_{r,n, r+b-n} 
\al
so that $\mathrm{H}_{n,r,b}$ is already determined by $\{n,r\}$ together with $r+b$.
Conversely and more precisely, we have:  

 {\rm(i)} $\mathrm{H}_{n,r,b}=\delta_a$ for some $a$ iff 
$n\wedge r\wedge b\wedge(r+b-n)=0$ and $n\wedge r =a$;

 {\rm(ii)}  $\mathrm{H}_{n,r,b}=\mathrm{B}_{1,p}$ for some $p\in\mathopen]0,1\mathclose]$
iff $n\wedge r=1$ and $\frac{n\vee r}{r+b}=p$;

 {\rm(iii)} in all other cases,  $\mathrm{H}_{n,r,b}$ is not a binomial law and determines
$\{n,r\}$ and $r+b$, that is, $\mathrm{H}_{n,r,b}=\mathrm{H}_{n',r',b'}$
for some $r',b'\in\N_0$ and $n'\in\{0,\ldots,r'+b'\}$ holds iff
$\{n,r\} = \{n',r'\}$ and $r+b=r'+b'$.

\smallskip\noindent{\rm\textbf{(c) Reflections.}}
$\mathrm{h}_{n,r,b}(k)=\mathrm{h}_{n,b,r}(n-k)$ for $k\in \Z$.

\smallskip\noindent{\rm\textbf{(d) Symmetries.}}
$\mathrm{H}_{n,r,b}$ is symmetric about its mean $\mu$ iff 
$n\wedge r\wedge b\wedge(r+b-n)=0$
or  $\frac{r+b}2\in\{n,r\}$, which  is the case  
iff $\kappa_3=0$, and which implies that $\mu \in\{\frac{n}2, \frac{r}2\}$.
\end{Lem}
\begin{proof} (a) Claim~\eqref{Eq:Supp_Hyp} is obvious from~\eqref{Eq:H_nrb_densitiy}.
The formulas for $\mu$ and $\sigma^2$ in~\eqref{Eq:Hyp_mean_var_kappa}
are proved in several textbooks as in \cite{Cornfield}, by considering 
a sum of indicator variables indicating ``red'' at each  of the $n$ draws, 
and this method works for $\kappa_3$ as well; alternatively one may 
use~\eqref{Eq:H_nrb_densitiy} and the 
differential equation for hypergeometric functions as in \cite[\S~5.14]{Kendall}. 

(b)  With  $\alpha \coloneqq n\wedge r$, $\beta\coloneqq n\vee r$, and $N\coloneqq r+b$,
a computation starting from~\eqref{Eq:H_nrb_densitiy} yields
\la               \label{Eq:Hyp_parameters_unidentifiable_2}
 \mathrm{h}_{n,r,b}(k) 
 &=& \frac{n_{}^{\underline{k}} r_{}^{\underline{k}} b_{}^{\underline{n-k}} }{k!(r+b)_{}^{\underline{n}} }
 \,\ =\,\
\frac{\alpha_{}^{\underline{k}}\beta_{}^{\underline{k}}
(N-\beta)_{}^{\underline{\alpha-k}}  }{k!N_{}^{\underline{\alpha}}} \quad\text{ for }k\in\{0,\ldots,n\},
\al
hence \eqref{Eq:Hyp_parameters_unidentifiable_1}.

(i) follows from \eqref{Eq:Supp_Hyp} and the formula for $\sigma^2$ 
in~\eqref{Eq:Hyp_mean_var_kappa}.

(ii) The ``if'' claim is clear by~\eqref{Eq:Hyp_parameters_unidentifiable_2}
with $k\in\{0,1\}$. Conversely, 
if $\mathrm{H}_{n,r,b}=\mathrm{B}_{1,p}$ with $p\in\mathopen]0,1\mathclose]$,
then $n\wedge r =1$ by~\eqref{Eq:Supp_Hyp}, and 
$p=\mu =\frac{n\vee r}{r+b}$ in view of~\eqref{Eq:Hyp_mean_var_kappa}.  
 
(iii) Assume that $\mathrm{H}_{n,r,b}$ is not as in (i) or (ii) and, 
without loss of generality in view of~\eqref{Eq:Hyp_parameters_unidentifiable_1},
that $n\le r$. Then $r\wedge b >0$ and $n>1$, hence also $0<\mu<n$,
and~\eqref{Eq:Hyp_mean_var_kappa}  yields 
\la                                                         \label{Eq:mu_sigma_supsupp}
  \sigma^2 &=& \mu\,\left(1-\frac{\mu}{n} \right)\frac{r+b-n}{r+b-1}
  \,\ <\,\ \mu\, \left( 1 - \frac{\mu}n\right). 
\al
The identiy in \eqref{Eq:mu_sigma_supsupp} yields $r+b$ as a function of the mean 
$\mu$, the variance $\sigma^2$, and  the right endpoint $n=n\wedge r$ 
of  $\mathrm{H}_{n,r,b}$, and then
$r = (r+b)\mu/n$ and hence $\{n,r\}$ as a function of quantities already determined
by   $\mathrm{H}_{n,r,b}$. The inequality 
$\sigma^2 < \mu\, \left( 1 - \frac{\mu}n\right)$, as a  relation 
between the mean, the variance, and the right endpoint of a law, would instead be an equality
if~$\mathrm{H}_{n,r,b}$ were binomial. 

(c) Trivial using~\eqref{Eq:H_nrb_densitiy}. 

(d) If $\mathrm{H}_{n,r,b}$ is symmetric about its mean, then $\kappa_3=0$,
as for any law with existing third moment.
If  $\kappa_3=0$, then~\eqref{Eq:Hyp_mean_var_kappa} yields
the stated condition for the parameters. If the latter holds,
then $\sigma^2=0$ and symmetry is trivial, 
or $\frac{r+b}2\in\{n,r\}$ and then (c) yields for $k\in\Z$ either 
$r=b$ and hence
\[
 \mathrm{h}_{n,r,b}(k) &=& \mathrm{h}_{n,b,r}(n-k) 
\,\ =\,\  \mathrm{h}_{n,r,b}(n-k),
\]
or $n=r+b-n$ and hence, using also~\eqref{Eq:Hyp_parameters_unidentifiable_1}
at the first and at the last step below, 
$\mathrm{h}_{n,r,b}(k) =\mathrm{h}_{r,n,r+b-n}(k)=\mathrm{h}_{r,r+b-n,n}(r-k)
 = \mathrm{h}_{r,n,r+b-n}(r-k) = \mathrm{h}_{n,r,b}(r-k)$, 
and hence in either case
the symmetry of $\mathrm{H}_{n,r,b}$, necessarily about its mean.
The final claim about $\mu$ is obvious using~\eqref{Eq:Hyp_mean_var_kappa}.
\end{proof}

Let $P$ be a binomial or a hypergeometric law. We then call $N\in\N\cup\{\infty\}$ 
a {\em population size parameter} of $P$ if $N=\infty$ and $P$ is binomial,
or if $P=\mathrm{H}_{n,r,b}$ for some $r,b\in\N_0$ and $n\in\{0,\ldots,r+b\}$
with $r+b=N$. 
By Lemma~\ref{Lem:Hyp_basics}(b), $N$ is uniquely determined by $P$ unless
$P$ is a Dirac or a Bernoulli law. Given a population size parameter 
$N$ of $P$, we let $\sigma_0^2$ denote the 
{\em usual approximate variance} of $P$, with respect to $N$, namely,
with $\sigma^2$ denoting the true variance of $P$,
\la                                                   \label{Eq:Def_usual_appr_var}  
 \sigma_0^2 
  &\coloneqq & \left\{\begin{array}{ll}  0 &\quad\text{ if }N=0,\\
    \frac{N-1}{N} \sigma^2 &\quad\text{ if }N\in\N, \\
    \sigma^2\phantom{\int\limits^1    }       &\quad \text{ if }N=\infty,  \end{array}\right. 
\al
which is uniquely determined by $P$, and hence may then be denoted by $\sigma_0^2(P)$,
 unless $P=\mathrm{B}_{1,p}$ with $p\in\mathopen]0,1\mathclose[$.
The customary but somewhat illogical dependence of $\sigma_0^2$ 
not only on $P$ in this last case 
is a source of the slightly awkward ``except'' proviso at the end of 
Theorem~\ref{Thm:Main}(a) below.

%%%%%%%%%%%%%%%%%%%%%%%%%%%%%%%%%%%%%%%%%%%%%%%%%%%%%%%%%%%%%%%%%%%%%%%%%%%%%%%%%%%%%%%%%%
%%%%%%%%%%%%%%%%%%%%%%%%%%%%%%%%%%%%%%%%%%%%%%%%%%%%%%%%%%%%%%%%%%%%%%%%%%%%%%%%%%%%%%%%%%%
\subsection{The main result} 
%%%%%%%%%%%%%%%%%%%%%%%%%%%%%%%%%%%%%%%%%%%%%%%%%%%%%%%%%%%%%%%%%%%%%%%%%%%%%%%%%%%%%%%%%%%
%%%%%%%%%%%%%%%%%%%%%%%%%%%%%%%%%%%%%%%%%%%%%%%%%%%%%%%%%%%%%%%%%%%%%%%%%%%%%%   Theorem 1.3
\begin{Thm}                                                                 \label{Thm:Main}
%%%%%%%%%%%%%%%%%%%%%%%%%%%%%%%%%%%%%%%%%%%%%%%%%%%%%%%%%%%%%%%%%%%%%%%%%%%%$%%%%%$$$$$%%%%%
{\rm\textbf{(a)}}
Let $F$ and $f$ be the distribution function and the density of a symmetric 
hypergeometric or symmetric binomial law, with mean $\frac{n}2$, standard deviation $\sigma>0$, 
population size parameter $N$, and the usual approximate standard deviation $\sigma_0$. 
Let $G$ be the distribution function of a normal law with mean  $\frac{n}2$ and standard 
deviation $\tau\in[\sigma_0,\sigma]$. Then, for $s\in\R$, 
\la                                            \label{Eq:Main_ineq}
 \left| F(s)-G(s)\right| \, <\, d \,
 \text{ if }s\neq\left\lfloor\textstyle\frac{n}2\right\rfloor
  &\text{ and }& 
 \left| F(s-)-G(s-)\right| \,<\,d  \,
 \text{ if }s\neq\left\lceil\textstyle\frac{n}2\right\rceil
\al
holds with 
\la                                           \label{Eq:Def_d}
 d &\coloneqq &  F\left(\left\lfloor\textstyle\frac{n}2\right\rfloor \right) 
       -G\left(\left\lfloor\textstyle\frac{n}2\right\rfloor \right)
   \,\ =\,\ 
    G\left({\left\lceil\textstyle\frac{n}2\right\rceil}- \right) 
   -F\left({\left\lceil\textstyle\frac{n}2\right\rceil}-\right) 
  \,\ =\,\ \left\|F-G\right\|_\infty   \\ 
  &=&                                           \label{Eq:d_n_odd_even}
    \left\{ 
      \begin{array}{ll}                                     
       \Phi\left(\frac1{2\tau}\right) - \frac12 & \text{ if $n$ is odd},\\
       \frac12 f\left(\frac{n}2 \right)  & \text{ if $n$ is even}
      \end{array}
    \right\} \\                                            \label{Eq:d_in_interval}
  &\in& \frac1{\sigma} \cdot 
    \left[ \frac{\Phi(1)-\frac12}{2}\,,\, \frac1{\sqrt{8\pi}} \right[             
  \,\  = \,\ \left[\frac{0.17\!\ldots}{\sigma} \,,\, \frac{0.19\ldots}\sigma \right[  
\al
except that the upper bound claim $d<\frac{1}{\sigma\sqrt{8\pi}}$ 
in \eqref{Eq:d_in_interval} is false if we have both
$N=2$ and $ \tau/\sigma %\frac{\tau}{\sigma}
\le c =0.78\ldots$, %c =0.783923\ldots
with~$c$ defined by 
$\sqrt{2\pi}\left(\Phi\left(\frac1{c}\right)-\frac12\right)=1$.

\smallskip\noindent{\rm\textbf{(b)}} 
The interval $\left[ \frac{\Phi(1)-\frac12}{2}\,,\, \frac1{\sqrt{8\pi}} \right[$
in part~{\rm(a)} is the least possible, even if we 
assume there in addition that $N=\infty$ (binomial case) 
and hence $\tau=\sigma=\sigma_0=\frac12\sqrt{n}$.   
\end{Thm}
%%%%%%%%%%%%%%%%%%%%%%%%%%%%%%%%%%%%%%%%%%%%%%%%%%%%%%%%%%%%%%%%%%%%%%%%%%%%%%%%%%%%%%%%%%
Theorem~\ref{Thm:Main} and the supplements stated in the following 
Remark~\ref{Rem:Complements_to_main_result} are proved at the end of this paper.
%%%%%%%%%%%%%%%%%%%%%%%%%%%%%%%%%%%%%%%%%%%%%%%%%%%%%%%%%%%%%%%%%%%%%%%%%%%%%%%%%%%%%%
\begin{Rem}                                  \label{Rem:Complements_to_main_result}
%%%%%%%%%%%%%%%%%%%%%%%%%%%%%%%%%%%%%%%%%%%%%%%%%%%%%%%%%%%%%%%%%%%%%%%%%%%%%%%%%%%%%%%
{\rm\textbf{(a)}} If $\tau$ is restricted to be $\sigma$ in Theorem~{\rm\ref{Thm:Main}},
then the formulation obviously simplifies a bit, and in particular 
the ``except'' proviso concerning~\eqref{Eq:d_in_interval}
becomes redundant. 

\smallskip\noindent{\rm\textbf{(b)}}
Under the assumptions of Theorem~{\rm\ref{Thm:Main}(a)} as stated,
and if $\frac{N}{N-1}$ is read as $1$ in case of $N=\infty$, 
we have without any exception 
\la                                                    \label{Eq:d-tau-bounds}
 \frac{\Phi\left(\sqrt{2}\right)-\frac12 }{\sqrt{8}\,\tau}   
  &\le&\frac{\sqrt{\frac{N-1}{N}}\left(\Phi\left(\sqrt{\frac{N}{N-1}}\right)
      -\frac12\right) }{2\tau}    
  \,\ \le \,\ d \,\ < \,\ \frac{1}{\tau\sqrt{8\pi}}.
\al
\noindent{\rm\textbf{(c)}} In the special case of $\tau =\sigma_0$, the upper bound for 
$d$ in~\eqref{Eq:d-tau-bounds} can be refined to 
\la                                  \label{Eq:d-sigma_null-bounds}
 d  \,\ \le\,\ \Phi\left(\frac1{2\sigma_0}\right) -\frac12
  \,\ < \,\ \frac{1}{\sqrt{8\pi}\sigma_0},
\al
with equality in the first inequality iff $n$ is odd.
\end{Rem}

Theorem~{\rm\ref{Thm:Main} specialized to symmetric binomial laws with $N=\infty$ 
reduces to \cite[Theorem 1.1 and Corollaries 1.1 and 1.2]{HM}. All other results
in the literature related to Theorem~\ref{Thm:Main} and known to us yield weaker 
or incomparable conclusions under more general hypotheses. Let us mention a few of these:

The central limit theorem for hypergeometric laws, namely 
``$\left\|F-G\right\|_\infty \rightarrow 0$  if $\sigma\rightarrow\infty$'' with the notation
of Theorem~\ref{Thm:Main} extended to not necessarily symmetric laws, is proved by 
R\'enyi in~\cite[pp.~ 465--466]{Renyi} as a corollary to~\cite{ErdosRenyi}.   
R\'enyi names S.N.~Bernstein 
as the originator under the  additional assumption
``$\frac{r}{r+b}$ constant'' in the notation of subsection~\ref{Lem:Hyp_basics}.
He also states that a direct proof of the general case ``leads to tiresome 
calculations'', which is  refuted by   
Morgenstern's treatment in \cite[pp.~62--63]{Morgenstern},
where the appropriate local central limit theorem is elegantly derived
from the corresponding one for binomial laws 
by writing $\mathrm{h}_{n,r,b}(k)=\mathrm{b}_{r,p}(k)\mathrm{b}_{b,p}(n-k)/\mathrm{b}_{r+b}(n)$
with $p\coloneqq \frac{n}{r+b}$ in the notation of subsection~\ref{Lem:Hyp_basics}.

Let now $C$ denote the optimal Berry-Esseen constant in the non-i.i.d.~case, 
so that $0.4097<(3+\sqrt{10})/(6\sqrt{2\pi})\le C  < 0.5583$ with the upper bound as announced 
in~\cite{Shev2013}. 
Let further $F$ be the distribution function of a 
Bernoulli convolution  $P = \bigconv_{j=1}^n\mathrm{B}_{p_j}$ with $p\in[0,1]^n$, and let 
$G$ be the distribution function of a normal law with the same mean $\mu=\sum_{j=1}^np_j$
and variance $\sigma^2= \sum_{j=1}^np_j\left(1-p_j\right)$.
Then, since $\beta_j:= p_j(1-p_j)( p_j^2+(1-p_j)^2)\le
 p_j\left(1-p_j\right) $ is the 
third absolute moment of $\mathrm{B}_{p_j}$, we have
$\left\| F -G \right\|_\infty  \le C\sigma^{-3}\sum_{j=1}^n \beta_j$  
and hence
\la         \label{BE_for_BCs} 
  \frac{1}{2\sqrt{1+12\sigma^2}} &\le & \left\| F -G \right\|_\infty 
   \,\ < \,\ \frac{0.5583}{\sigma},
\al
where the lower bound follows from the continuity of $G$ 
and from the lower bound for the maximal jump size of $F$  
obtained from~\eqref{Eq:Conc_vs_sigma} below with $h=1$.
Now it is well known from \cite[Corollary~5 with $n=2$, 
hence $F_2$ generating function of $\mathrm{H}_{s_1,s_2,N-s_2}= \mathrm{H}_{s_2,s_1,N-s_1}$]{Vat_Mik}
that every hypergeometric law is a Bernoulli convolution as above, 
with certain in general not explicitly available $p_j$, but of course $\mu$
and $\sigma^2$  computable from~\eqref{Eq:Hyp_mean_var_kappa}.
Thus, as already known from~\cite[Theorem 1 with $n=2$, rewritten in terms 
of $\mu_0+s_1+s_2-N$]{Vat_Mik} in case of the upper bound,
\eqref{BE_for_BCs} directly applies to $F$ and $G$ as in the previous paragraph, 
and thus yields a result more explicit than the two theorems in~\cite{LC}
and with a simpler proof, but~\eqref{BE_for_BCs}  is in the symmetric case of 
course weaker than~\eqref{Eq:d_in_interval} applied to $\left\| F -G \right\|_\infty$. 

H\"oglund's theorem already mentioned in subsection~\ref{Subsec:Background}
yields the upper bound in~\eqref{BE_for_BCs},  in the general hypergeometric case,
with an unspecified constant in place of $0.5583$. 

Some further related results and references
can be found in the papers
\cite{Mohamed_Mirakhmedov} concerning in particular sums of simple random samples,
\cite{Neammanee} concerning Bernoulli convolutions,
and~\cite{LCM} concerning hypergeoemetric laws.

%%%%%%%%%%%%%%%%%%%%%%%%%%%%%%%%%%%%%%%%%%%%%%%%%%%%%%%%%%%%%%%%%%%%%%%%%%%%%%%%%%%%%%%%%%%%%%%
\subsection{On concentration-variance inequalities}
%%%%%%%%%%%%%%%%%%%%%%%%%%%%%%%%%%%%%%%%%%%%%%%%%%%%%%%%%%%%%%%%%%%%%%%%%%%%%%%%%%%%%%%%%%%%%%%
In deriving the lower bound in~\eqref{BE_for_BCs} 
above, we have used inequality~\eqref{Eq:Conc_vs_sigma} below, which is   
due to Paul L\'evy in a sharper version.
%%%%%%%%%%%%%%%%%%%%%%%%%%%%%%%%%%%%%%%%%%%%%%%%%%%%%%%%%%%%%%%%%%%%%%%%%%%%%%%%%%%%%%%%%%% 
\begin{Lem}                                                         \label{Lem:Conc_vs_sigma} 
%%%%%%%%%%%%%%%%%%%%%%%%%%%%%%%%%%%%%%%%%%%%%%%%%%%%%%%%%%%%%%%%%%%%%%%%%%%%%%%%%%%%%%%%%%%
Let $P$ be a law on $\R$ with variance $\sigma^2$. Then we have
\la
 \sup_{x\in\R}P(\mathopen]x,x+h\mathclose[) &\ge& \frac{h}{\sqrt{h^2+12\sigma^2}}\quad\text{ for }
       h\in\mathopen]0,\infty\mathclose[,     \label{Eq:Conc_vs_sigma} \\ 
 \underset{x\in\R}{\mathrm{ess\,sup\,}} f(x)                           \label{Eq:esssup_vs_sigma}
 &\ge & \frac{1}{\sqrt{12\sigma^2}} \quad\text{ if $f$ is a Lebesgue density of $P$}.
\al
\end{Lem}
\begin{proof} For~\eqref{Eq:esssup_vs_sigma} we may assume that $P$ has mean zero
and $M:=\text{L.H.S.\eqref{Eq:esssup_vs_sigma}}<\infty$. 
With $c\coloneqq \frac1{2M}$, we then  have $\alpha \coloneqq \int_{|x|>c}f(x)\dd{x}
=\int_{|x|\le c}(M-f(x))\dd{x}$, hence 
$\int_{|x|>c}x^2f(x)\,\mathrm{d}x\ge c^2\alpha\ge \int_{|x|\le c}x^2\left(M-f(x)\right)\mathrm{d}x$, 
and thus
$\sigma^2=\int x^2f(x)\dd{x}\ge \int_{|x|\le c}x^2\left(f(x)+M-f(x)\right)\mathrm{d}x$
$=$ $\frac23c^3M=\frac{1}{12M^2}$.
 
To prove now~\eqref{Eq:Conc_vs_sigma}, we apply~\eqref{Eq:esssup_vs_sigma}
to the density
$x\mapsto g(x) \coloneqq\frac1h P(\mathopen]x-h,x\mathclose[)$ 
and  the variance $\sigma^2+h^2/12$
of the convolution of $P$ with the uniform law on $]0,h[$, 
to get  
% \la                  \label{Eq:Proof_conc_inq}
%  \text{L.H.S.\eqref{Eq:Conc_vs_sigma}} &=& 
%   h\sup_{x\in\R} g(x) \,\ \ge \,\ h \,\underset{x\in\R}{\mathrm{ess\,sup\,}} g(x)
%   \,\ \ge  \,\ h  \frac{1}{ \sqrt{ 12\left( \sigma^2 + \frac{h^2}{12}  \right)    }  }. 
% \al
\begin{equation}                  \label{Eq:Proof_conc_inq}
 \text{L.H.S.\eqref{Eq:Conc_vs_sigma}} \,\ = \,\  
  h\sup_{x\in\R} g(x) \,\ \ge \,\ h \,\underset{x\in\R}{\mathrm{ess\,sup\,}} g(x)
  \,\ \ge  \,\ h  \frac{1}{ \sqrt{ 12\left( \sigma^2 + \frac{h^2}{12}  \right)    }  }. \qedhere
\end{equation}
\end{proof}

L\'evy \cite[p.~149, Lemme 48,1]{Levy} proved under the assumption of Lemma~\ref{Lem:Conc_vs_sigma}:
If $p\in\N$  and  $\lambda \in [0,1]$ are such that 
$c\coloneqq \text{L.H.S.\eqref{Eq:Conc_vs_sigma}}=\frac{\lambda}p +\frac{1-\lambda}{p+1}$, then 
\la                                                       \label{Eq:Conc_vs_sigma_Levy}
  12\frac{\sigma^2}{h^2} &\ge& \lambda p^2 +(1-\lambda)(p+1)^2 - 1,
\al
with equality for 
$P = \frac{\lambda}{p}\sum_{j=0}^{p-1}\delta_{jh} +  \frac{1-\lambda}{p+1}\sum_{j=0}^{p}\delta_{(j-\frac12)h}$.
Writing $p^2=\left(\frac1p \right)^{-2}$ and $(p+1)^2=\left(\frac1{p+1} \right)^{-2}$, 
and using convexity, \eqref{Eq:Conc_vs_sigma_Levy} yields 
$ 12\frac{\sigma^2}{h^2} \ge \left(\frac{\lambda}p + \frac{1-\lambda}{p+1} \right)^{-2}-1 = c^{-2} -1$,
hence~\eqref{Eq:Conc_vs_sigma}, and~\eqref{Eq:esssup_vs_sigma}
follows easily using  $\frac1hP(\mathopen]x,x+h\mathclose[) \le \mathrm{ess\,sup\,}f$. 
We refer to \cite[p.~27]{Hen_Theo} for a proof of~\eqref{Eq:Conc_vs_sigma_Levy}
more formal than L\'evy's,   and to \cite{Foley_Hill_Spruill} for generalizations.

The present proof of first~\eqref{Eq:esssup_vs_sigma} and then~\eqref{Eq:Conc_vs_sigma}  
is a slightly simplified and corrected version of an argument given 
by Bobkov and Chistyakov: Our first part is simpler, or at least more elementary,
than~\cite[first 5 lines of Proof of Proposition~2.1]{Bobkov_Chist}.
To see the correction in the second part, let us first observe that 
we actually have equality at the second step in~\eqref{Eq:Proof_conc_inq}, 
since our $g$ is lower semicontinuous, but that  this could be wrong if we had  closed intervals 
$[x,x+h]$ on the left in~\eqref{Eq:Conc_vs_sigma} and analogously 
also in the  definition of $g$, as for example if $P=\frac12(\delta_0+\delta_1)$ 
and $h=1$, contrary to~\cite[(2.1)]{Bobkov_Chist} where hence $Q(X;\lambda)$
should be replaced by $Q(X;\lambda{-})$.

Finally we have to mention
that  \eqref{Eq:esssup_vs_sigma} also follows by letting $p\rightarrow \infty$
in Moriguti's sharp inequality~\cite[(3.4)]{Moriguti} for $\mathrm{L}^p$-norms,
valid under the hypothesis of~\eqref{Eq:esssup_vs_sigma}, namely 
\[
  \left\|f\right\|_p &\ge& 
   \left( \tfrac{2p}{3p-1} \right)^{\frac1p}  
   \left(\sqrt{\tfrac{p-1}{3p-1}}/\left(\mathrm{B}\left(
    \tfrac{p}{p-1},\tfrac12 \right)\sigma \right) \right)^{1-\frac1p} 
  \quad \text{ for }p\in\mathopen]1,\infty\mathclose[.
\]

%%%%%%%%%%%%%%%%%%%%%%%%%%%%%%%%%%%%%%%%%%%%%%%%%%%%%%%%%%%%%%%%%%%%%%%%%%%%%%%%%%%%%
\subsection{The method of proof}
%%%%%%%%%%%%%%%%%%%%%%%%%%%%%%%%%%%%%%%%%%%%%%%%%%%%%%%%%%%%%%%%%%%%%%%%%%%%%%%%%%%%%
The proof of Theorem~\ref{Thm:Main} near the end of  section~\ref{Sec:4} below rests on the 
following simple lemma, which was implicitly used  also in \cite{HM}.
%%%%%%%%%%%%%%%%%%%%%%%%%%%%%%%%%%%%%%%%%%%%%%%%%%%%%%%%%%%%%%%%%%%%%%%%%%
\begin{Lem}                 \label{Lem:Symmetry_and_distances}           %  Lemma 1.3
%%%%%%%%%%%%%%%%%%%%%%%%%%%%%%%%%%%%%%%%%%%%%%%%%%%%%%%%%%%%%%%%%%%%%%%%%%
Let $F$ and $G$ be distribution functions of laws $P$ and $Q$ on $\R$ 
with $P(\Z)=1$, $G$ continuous and strictly increasing,  $P$ and $Q$ symmetric
about $\frac{n}2\in\R$, and $d$ defined by the first equality in~\eqref{Eq:Def_d}.
Then we have $n\in \Z$, the second equality in~\eqref{Eq:Def_d}, and 
\la                                                       \label{Eq:Main_ineq_0}
  d &=& \left\{\begin{array}{ll} G\left(\frac{n+1}2\right)-\frac12 &\text{ if $n$ is odd,}\\
          \frac12 P \left(\left\{\frac{n}2\right\}\right) & \text{ if $n$ is even.}
        \end{array}\right.  
\al
Further, \eqref{Eq:Main_ineq} holds for every $s\in\R$ iff the following two conditions
are satisfied:
\la
 F(s)-G(s)                                                           \label{Eq:Main_ineq_1}
    &<&d\quad\text{ for $s\in\Z$ with $s>\left\lfloor\textstyle\frac{n}2\right\rfloor$},\\
 G(s)-F(s-1)                                                         \label{Eq:Main_ineq_2}
    &<&d\quad\text{ for $s\in\Z$ with $s>\left\lceil\textstyle\frac{n}2\right\rceil$}.
\al 
\end{Lem}
\begin{proof} The symmetry assumptions can be written as
\la                                                                  \label{Eq:P_Q_sym}
  F(s)\,=\,1-F((n-s)-)&\text{ and }& G(s)\,=\,1-G((n-s)-) \quad\text{ for }s\in\R.
\al
The assumption $P(\Z)=1$ then yields $0<P(\{k\})=F(k)-F(k-1)=P(\{n-k\})$ for some 
$k\in\Z$, and hence $n\in\Z$. Next, \eqref{Eq:P_Q_sym} for 
$s=\left\lfloor\textstyle\frac{n}2\right\rfloor$ yields 
\[
     F\left(\left\lfloor\textstyle\frac{n}2\right\rfloor\right) 
  +  F\left(\left\lceil\textstyle\frac{n}2\right\rceil-\right)
 &=& 1 \,\ =\,\ G\left(\left\lfloor\textstyle\frac{n}2\right\rfloor\right) 
  +  G\left(\left\lceil\textstyle\frac{n}2\right\rceil-\right)
\]
and we get the second equality in~\eqref{Eq:Def_d}, and also 
$G\left(\left\lceil\textstyle\frac{n}2\right\rceil-\right)
- F\left(\left\lceil\textstyle\frac{n}2\right\rceil-\right)=G\left(\frac{n+1}2\right)-\frac12$
if $n$ is odd, and 
$F\left(\left\lfloor\textstyle\frac{n}2\right\rfloor\right) 
- G\left(\left\lfloor\textstyle\frac{n}2\right\rfloor\right) = 
\frac12 + \frac12P\left(\left\{\frac{n}2\right\} \right) -\frac12$ if $n$ is even,
and thus~\eqref{Eq:Main_ineq_0}. 

Trivially, \eqref{Eq:Main_ineq} implies~\eqref{Eq:Main_ineq_1} and~\eqref{Eq:Main_ineq_2}.
Conversely, let us assume~\eqref{Eq:Main_ineq_1} and~\eqref{Eq:Main_ineq_2}. 
If $s\in\R\setminus\Z$, then 
$F(s)-G(s)=F(\lfloor s\rfloor)-G(s) <F(\lfloor s\rfloor)-G(\lfloor s\rfloor)$
and  $G(s)-F(s)=G(s)-F(\lceil s\rceil-)< G(\lceil s\rceil- )-F(\lceil s\rceil-)$; hence it is
enough to prove~\eqref{Eq:Main_ineq} for $s\in\Z$. 
If $s\in\Z$ with $s>\left\lfloor\textstyle\frac{n}2\right\rfloor$, 
then $F(s)-G(s)<d$ by~\eqref{Eq:Main_ineq_1}, and   $G(s)-F(s)
\le G(s+1)-F(s)<d$ by~\eqref{Eq:Main_ineq_2} as $s+1>\left\lceil\textstyle\frac{n}2\right\rceil$;
hence $|F(s)-G(s)|<d$. If $s\in\Z$ with $s<\left\lfloor\textstyle\frac{n}2\right\rfloor$, 
then $t\coloneqq n-s>\left\lceil\textstyle\frac{n}2\right\rceil$, 
and~\eqref{Eq:P_Q_sym}, \eqref{Eq:Main_ineq_1}, \eqref{Eq:Main_ineq_2} yield
$F(s)-G(s)= G(t)-F(t-1)<d$ and $G(s)-F(s)= F(t-1)-G(t)\le F(t)-G(t)<d$;
hence again $|F(s)-G(s)|<d$. Thus the first part of~\eqref{Eq:Main_ineq} holds for $s\in\Z$, 
and the second follows by applying, for a given 
$s\neq \left\lceil\textstyle\frac{n}2\right\rceil$, 
the first one to $t \coloneqq  n-s \neq \left\lfloor\textstyle\frac{n}2\right\rfloor$. 
\end{proof}

In the  situation of Theorem~\ref{Thm:Main}, assumption~\eqref{Eq:Main_ineq_1}
and part of assumption~\eqref{Eq:Main_ineq_2} are proved below in 
Lemmas~\ref{Lem:4.2} and~\ref{Lem:4.3} by monotonicity considerations, 
and the part of~\eqref{Eq:Main_ineq_2} not thus covered is proved 
by using lower bounds for $d$ from  Lemma~\ref{Lem:Bounds_h_etc_n_even_Neu}  
together with Lemma~\ref{Lem:4.4}.
The proofs of the lemmas of section~\ref{Sec:4}
use various auxiliary inequalities from sections~\ref{Sec:Aux_inequalities}
and~\ref{Sec:Normal_df_inequalities}.

%%%%%%%%%%%%%%%%%%%%%%%%%%%%%%%%%%%%%%%%%%%%%%%%%%%%%%%%%%%%%%%%%%%%%%%%%%%%%%%%%%%%%%
\section{Some standard analytic inequalities}                  \label{Sec:Aux_inequalities}
%%%%%%%%%%%%%%%%%%%%%%%%%%%%%%%%%%%%%%%%%%%%%%%%%%%%%%%%%%%%%%%%%%%%%%%%%%%%%%%%%%%%%%
Very elementary inequalities, like $1+x<\mathrm{e}^x$ for $x\in\R\setminus\{0\}$  and 
$\frac{x}{1+x}<\log(1+ x)<x$ for $x>-1$,  will often be used without comment. 
%%%%%%%%%%%%%%%%%%%%%%%%%%%%%%%%%%%%%%%%%%%%%%%%%%%%%%%%%%%%%%%%%%%%%%%%   Lemma 2.1  
\begin{Lem}                                                  \label{Lem:(1+x)/(1+y)}
%%%%%%%%%%%%%%%%%%%%%%%%%%%%%%%%%%%%%%%%%%%%%%%%%%%%%%%%%%%%%%%%%%%%%%%%%%%%%%%%%%%%%%             
If $x,y\in\R$ satisfy $0\le y\le|x|$ or $x\le y\le 0$ 
or $x-\frac23x^2\ge -y\ge 0$, then
\la
  (1+x)\mathrm{e}^{-x} &\le& (1+y)\mathrm{e}^{-y},
\al
and equality holds iff $x=y$. The constant $\frac23$ in the 
assumption can not be lowered.
\end{Lem}
\begin{proof} See  \cite[Lemma~2.1]{HM}.
\end{proof}
%%%%%%%%%%%%%%%%%%%%%%%%%%%%%%%%%%%%%%%%%%%%%%%%%%%%%%%%%%%%%%%%%%%%%%
\begin{Lem}                                                           \label{Lem:cosh}
%%%%%%%%%%%%%%%%%%%%%%%%%%%%%%%%%%%%%%%%%%%%%%%%%%%%%%%%%%%%%%%%%%%%%%
Let %$0\neq x\in\R$.
$x\in\R\setminus\{0\}$. 
Then $\exp(\frac{x^2}2-\frac{x^4}{12})<\cosh(x)<(1+\frac{x^2}3)\exp(\frac{x^2}6)$.
\end{Lem}
\begin{proof} 
Analogously to \cite[Erster Abschnitt, Aufgabe 154 und L\"osung, pp.~28, 183]{PolyaSzegoeI}, 
the partial fraction expansion of the hyperbolic tangent function
\[
 \frac{\tanh(x)}{x} &=& \sum_{k=1}^\infty \frac2{((k+\frac12)\pi)^2+x^2}
\]
proved for example in~\cite[pp.~199, 294]{RemmertSchumacher} implies that
$\tanh(x)/x$ is enveloped by its power series around zero, namely
\[
(-)^{n}\left(\frac{\tanh(x)}{x}-\sum_{k=0}^n (-)^k\alpha_kx^{2k}  \right)
&>& 0 \quad\text{ for }x\in\R\setminus\{0\}, n \in \N_0,
\]
where $\alpha_0=1, \alpha_1=\frac13, \alpha_2=\frac2{15},\ldots$,
and using $\log(\cosh(x))=\int_0^x\tanh(t)\dd t$ then yields
\la        \label{Eq:logcosh_env}
 \quad (-)^{n}\left(\log(\cosh(x))- \sum_{k=0}^n (-)^k\frac{\alpha_k}{2k+2}x^{2k+2}\right)
 &<&0 \quad\text{ for } x\in\R\setminus\{0\}, n \in \N_0.
\al
Taking $n=1$ yields the first inequality claimed. 

To prove the second one, which improves the case $n=0$ of \eqref{Eq:logcosh_env},
we observe that
the coefficients of $x^{2k}$ in the power series of the two functions involved,
namely $a_k\coloneqq \frac1{(2k)!}$ % for $k\in\N_0$, $b_0\coloneqq 1$,  
and $b_k\coloneqq \frac{2k+1}{6^kk!}$ for $k\in\N_0$, are all $>$ $0$, %positive,
and their 
quotients $c_k\coloneqq a_k/b_k$ satisfy $c_0=c_1=1$
and $c_{k+1}/c_k = 3/(2k+3)<1 $ for $k\in\N$.   
\end{proof}
%%%%%%%%%%%%%%%%%%%%%%%%%%%%%%%%%%%%%%%%%%%%%%%%%%%%%%%%%%%%%%%%%%%%%%%%%%% Lemma 2.3
\begin{Lem}                                                           \label{Lem:w(x)}
%%%%%%%%%%%%%%%%%%%%%%%%%%%%%%%%%%%%%%%%%%%%%%%%%%%%%%%%%%%%%%%%%%%%%%%%%%%%%%%%%%%%%%
With $w(x)\coloneqq \frac{\Gamma(2x+1)}{\Gamma^2(x+1)}2^{-2x}
=\frac{\Gamma(x+\frac12)}{\sqrt{\pi}\Gamma(x+1)}$, we have 
\la                                                     \label{Eq:w_bounds}
 -\frac1{8x}\,\,\,<\,\,\, \log\left(\sqrt{\pi x} w(x)\right)
 &<&  -\frac1{8x} + \frac1{192x^3} 
 \qquad\text{ for }x\in\mathopen]0,\infty\mathclose[  \\ 
&\le& -\frac{23}{192x}\,\, \,\,\, \qquad\qquad\text{ for }x\in[1,\infty[. \label{Eq:w_for_x_ge_1}
\al
\end{Lem}
\begin{proof}[Two proofs] 
Inequality \eqref{Eq:w_for_x_ge_1} is of course trivial in view of $\frac1{x^3}\le\frac 1x$.   
Concerning~\eqref{Eq:w_bounds}:

For integer $x$, and only this case   will be needed in this 
paper,   \eqref{Eq:w_bounds} is proved by Everett
in~\cite[(10), with $W_n$ there being the present  $(\sqrt{\pi n}w(n))^2$]{Everett}    .
 
For general $x$, Sasv\'ari~\cite{Sasvari} presents the inequalities in~\eqref{Eq:w_bounds} 
as special cases of a more general corollary to a theorem yielding the monotonicity 
in $x$ of the error of each of the asymptotic expansions 
$\sum_{j=1}^N c_{r,j}x^{1-2j}$  of 
$\log\left(\sqrt{2\pi\frac{r-1}r} \left(\frac{(r-1)^{r-1}}{r^{r}}\right)^x 
 \binom{rx}{x}\right)$
for $\mathopen]0,\infty\mathclose[ \ni x \rightarrow \infty$,
with $ r\in\mathopen]1,\infty\mathclose[$ and $N\in\N_0$ fixed and here
$\binom{rx}{x} \coloneqq  \Gamma(rx+1)/\left(\Gamma(x+1)\Gamma((r-1)x)\right)$. 
Sasv\'ari's proof is short and elegant but, to get just~\eqref{Eq:w_bounds} 
and its analogues in Sasv\'ari's corollary, can even be shortened
a bit by using in his formula~(2) and in his notation
just ``$Q_s<0$'' rather than ``$Q_s$ increasing''.   
\end{proof}

Although not needed here, let us remark that numerical calculations
suggest that we have in fact 
$\sup_{x\in\mathopen[1,\infty\mathclose[} x\log\left(\sqrt{\pi x} w(x)\right)
= \log(\frac12\sqrt{\pi})= -0.1207\ldots < -\frac{23}{192} = -0.1197\ldots$.

%%%%%%%%%%%%%%%%%%%%%%%%%%%%%%%%%%%%%%%%%%%%%%%%%%%%%%%%%%%%%%%%%%%%%%%%%%%%
%%%%%%%%%%%%%%%%%%%%%%%%%%%%%%%%%%%%%%%%%%%%%%%%%%%%%%%%%%%%%%%%%%%%%%%%%%%%
\section{Normal distribution function inequalities}            \label{Sec:Normal_df_inequalities}
%%%%%%%%%%%%%%%%%%%%%%%%%%%%%%%%%%%%%%%%%%%%%%%%%%%%%%%%%%%%%%%%%%%%%%%%%%%%
%%%%%%%%%%%%%%%%%%%%%%%%%%%%%%%%%%%%%%%%%%%%%%%%%%%%%%%%%%%%%%%%%%%%%%%%%%%%
For comparing normal distribution function increments with their midpoint derivative 
approximations, we will need the rather sharp inequalities~\eqref{Eq:Incr_Phi_vs_phi} 
below, which improve the ones in~\cite[p.~322, Lemma~1]{Feller} and 
in~\cite[pp.~475--476, Lemma~1]{Nicholson} in an optimal way.
%%%%%%%%%%%%%%%%%%%%%%%%%%%%%%%%%%%%%%%%%%%%%%%%%%%%%%%%%%%%%%%%%%%%%%%%%%%%%%%%%%% Lemma 3.1
\begin{Lem}                                                       \label{Lem:Incr_Phi_vs_phi}       
%%%%%%%%%%%%%%%%%%%%%%%%%%%%%%%%%%%%%%%%%%%%%%%%%%%%%%%%%%%%%%%%%%%%%%%%%%%%%%%%%%%%%%%%%%%%
For  $x,h\in\R$ with $h\neq0$, we have
\la                                     \label{Eq:Incr_Phi_vs_phi}
 \exp \left(\tfrac{(x^2-1)h^2}{24}-\tfrac{x^4h^4}{960} \right) 
 &<& \frac{\Phi(x+\frac h2)-\Phi(x-\frac h2)}{h \phi(x) } 
  \,\ <\,\ \exp \left(\tfrac{(x^2-1)h^2}{24}+\tfrac{h^4}{1440} \right),
\al
and these inequalities are optimal for small $h$ 
in the sense that we have
\la                               \label{Eq:Incr_Phi_vs_phi_opt}
 \log\left(\tfrac{\Phi\left(x+\frac h2\right)-\Phi\left(x-\frac h2\right)}
                  {h \phi(x) } \right)
 &=& \tfrac{(x^2-1)h^2}{24} + \tfrac{(-x^4-4x^2+2)h^4}{2880}  + O(h^6)
  \quad\text{ for $x,h$ bounded,}
\al
with  $\max_{x\in\R}(-x^4-4x^2 +2)=2$ and $\min_{x\in\R}(-x^4-4x^2 +2)/x^4=-3$.
\end{Lem}
\begin{proof}  For $x,y\in\R$, let 
\[
 \epsilon_1(x,y) \,\ \coloneqq  \,\ -\frac{x^4y^4}{60} &\text{ and }&
 \epsilon_2(x,y) \,\ \coloneqq  \,\ \frac{y^4}{90}
\]
and, for $i\in\{1,2\}$, 
\[
 f_i(x,y) &\coloneqq & \textstyle
  \frac{\sqrt{2\pi}}2\left(\Phi(x+y)-\Phi(x-y) \right)
   - y\exp\left(\textstyle-\frac{x^2}2+\frac{(x^2-1)y^2}6+\epsilon_i(x,y) \right). 
\]
Noting  that~\eqref{Eq:Incr_Phi_vs_phi} is unaffected 
by sign changes of $x$ or $h$, and writing $y$ in place $h/2$, we 
have to prove for $x\ge 0$ and $y>0$ the inequalities 
\la \label{Eq:f1>0>f2}
 f_1(x,y) &>&0 \,\ >\,\ f_2(x,y). 
\al
Now $f_i(x,0)=0$  and, with a subscript $y$ denoting the  
partial derivative with respect to that variable,
\[
 \frac{f_{i,y}(x,y)}
{\exp\left(\frac{x^2y^2}{6} -\frac{x^2+y^2}2\right)}  
 &=& 
  \frac{\cosh(xy)}{\exp\left(\frac{x^2y^2}{6}\right)    } 
-{\textstyle \left(1 -\frac{y^2}3+y\epsilon_{i,y}(x,y) +\frac{x^2y^2}3
   \right)\exp\left(\frac{y^2}3 +\epsilon_i(x,y) \right) } \\
 &=:&  g^{}_{i}(x,y).
\]

For  $i=1$, we use the first inequality in Lemma~\ref{Lem:cosh} 
and $1+t<\mathrm{e}^t$ for $0\neq t\in\R$
to get 
 \[g^{}_{1}(x,y)&>& \exp\left(\frac{x^2y^2}3-\frac{x^4y^4}{12} \right) 
- \exp\left(  y\epsilon_{1,y}(x,y)  +\frac{x^2y^2}3 +\epsilon_1(x,y)  \right)
 \,\ = \,\ 0, 
\]
considering the cases $x\neq 0$ and $x=0$ separately to check the strict inequality,
and hence the first half of~\eqref{Eq:f1>0>f2}.

For  $i=2$, the second inequality in Lemma~\ref{Lem:cosh} 
and  $\frac{x^2y^2}3\exp(\frac{y^2}3 +\frac{y^4}{90})\ge\frac{x^2y^2}3$
yield
\[
 g^{}_{2}(x,y) &\le&   \textstyle 1 + \frac{x^2y^2}3 
  -\left(1 -\frac{y^2}3+ \frac{2y^4}{45}
   \right)\exp\left(\frac{y^2}3 +\frac{y^4}{90}  \right)-  \frac{x^2y^2}3 
 \,\ =\,\ 1- \exp(g(y^2)) 
\]
where, for $t\in\R$, 
\[ 
 g(t) &\coloneqq & \textstyle \frac{t}3+\frac{t^2}{90} +
\log\left(1-\frac {t}3 +\frac {2t^2}{45}\right)
\]
is well-defined with $g(0)=0$ and, for  $t>0$, satisfies
\[
 g'(t)&=& \textstyle \frac 13 +\frac t{45} +
 \frac {-\frac 13 +\frac {4t}{45}}{1-\frac {t}3 +\frac {2t^2}{45}}
 \,\,\,=\,\,\,  \frac {\frac {t^2}{135}+\frac {2t^3}{2025}}
 {1-\frac {t}3 +\frac {2t^2}{45}}  \,\,\,>\,\,\, 0
\]
and hence $g(t)>0$, yielding  $g^{}_{2}(x,y)<0$ and hence the second half of~\eqref{Eq:f1>0>f2}.

With the Hermite polynomials $\mathrm{H}_n$ given by  
$\mathrm{H}_n(x)=(-)^n\mathrm{e}^{x^2/2}\partial_x^{n}\mathrm{e}^{-x^2/2}$, 
in particular $\mathrm{H}_0(x)=1, \mathrm{H}_2(x) = x^2-1, \mathrm{H}_4(x)=x^4-6x^2+3$,
a Taylor expansion around $h=0$ shows that, for $x,h$ bounded, we have 
\[
 \frac{\Phi(x+\frac h2)-\Phi(x-\frac h2)}{h \phi(x) } 
 &=& 
 \sum_{j=0}^\infty \frac{\mathrm{H}_{2j}(x) h^{2j}}{2^{2j}(2j+1)!} 
 \,\ = \,\ 
  1 +\mathrm{H}_2(x)\frac{h^2}{24} +\mathrm{H}_4(x)\frac{h^4}{1920} + O(h^6)
   \label{Increments of Phi vs.phi,Power Series, first few}   
\]
and hence an application of $\log(1+y)=y-y^2/2 +O(y^3)$ for $y$ near zero 
and a short computation yield~\eqref{Eq:Incr_Phi_vs_phi_opt}.
\end{proof}

%%%%%%%%%%%%%%%%%%%%%%%%%%%%%%%%%%%%%%%%%%%%%%%%%%%%%%%%%%%%%%%%%%%%%%%%%%    Lemma 3.2
\begin{Lem}                                                  \label{Lem:Phi_near_zero}           
%%%%%%%%%%%%%%%%%%%%%%%%%%%%%%%%%%%%%%%%%%%%%%%%%%%%%%%%%%%%%%%%%%%%%%%%%%%%%%%%%%%%%%%%
Let $x>0$. Then
\la                                              \label{Eq:Phi_near_zero}
 x\mathrm{e}^{-\frac{x^2}6}&<&  \sqrt{2\pi}\left(\Phi(x)-{\textstyle\frac12}\right)
 \,\ <\,\  x\mathrm{e}^{-\frac{x^2}6+\frac{x^4}{90}} .
\al
\end{Lem}
\begin{proof} The claim results  if we apply~\eqref{Eq:Incr_Phi_vs_phi} to  
$(0,2x)$ in place of $(x,h)$.
\end{proof}

%%%%%%%%%%%%%%%%%%%%%%%%%%%%%%%%%%%%%%%%%%%%%%%%%%%%%%%%%%%%%%%%%%%%%%%%%%%%%
The following lemma often improves on~\cite[Lemma~2.2]{HM},  which 
yields~\eqref{Eq:Quot_Phi_incr_bounds} with $\alpha=1$ but with the upper bound replaced 
by $\exp\left(-\frac{y^2-x^2}{2}+\frac{|h|}{2}\left(|y|-|x|\right)\right)$,
and it always improves on~\cite{SasvariLindsey}, where~\eqref{Eq:Quot_Phi_incr_bounds}
with $\alpha=1$ is only obtained for $h=y-x$ 
and with $\beta=1-\frac{h^2}2 <1-\frac{h^2}{12}$.

%%%%%%%%%%%%%%%%%%%%%%%%%%%%%%%%%%%%%%%%%%%%%%%%%%%%%%%%%%%%%%%%%%%%%%%%%%  Lemma 3.3
\begin{Lem}                                          \label{Lem:Quot_Phi_incr_bounds}                  
%%%%%%%%%%%%%%%%%%%%%%%%%%%%%%%%%%%%%%%%%%%%%%%%%%%%%%%%%%%%%%%%%%%%%%%%%%%%%%%%%%%%%
Let $h\in\R\setminus\{0\}$. Then   
\la                      \label{Eq:Quot_Phi_incr_bounds}
\,\, {\textstyle\exp\left(-\alpha\tfrac{y^2-x^2}{2}\right)} 
&<& \frac{\Phi\big(y+\frac h2\big)-\Phi\big(y-\frac h2\big)}
         {\Phi\big(x+\frac h2\big)-\Phi\big(x-\frac h2\big)}
\,\ <\,\ {\textstyle \exp(-\beta\tfrac{y^2-x^2}{2})}
\quad\text{ if }\, |x|<|y| 
\al
holds with the optimal constants 
\la                                                   \label{Eq:Opt_alpha_beta}
  \alpha\,\coloneqq \,1&\text{ and }&  
  \beta \,\coloneqq \, \frac{\frac h2\exp\left(-\frac{h^2}8 \right)}
                    {\sqrt{2\pi}\left(\Phi(\frac{h}2)-{\frac12}\right)} 
        \,>\, \exp\left(-\tfrac{h^2}{12}-\tfrac{h^4}{1440}\right)  \,>\, 1-\tfrac{h^2}{12}.    
\al
\end{Lem}
\begin{proof} Since~\eqref{Eq:Quot_Phi_incr_bounds} and \eqref{Eq:Opt_alpha_beta}
are unaffected by sign changes of $x$ or $y$ or $h$, we may and do 
always assume that $0\le x<y$ and $h>0$ in this proof. For $\gamma\in\R$, let 
\[
 f_\gamma(x) &\coloneqq & \mathrm{e}^{\gamma\tfrac{x^2}2} 
                \left(\Phi\big(x+\tfrac h2\big)-\Phi\big(x-\tfrac h2\big)\right)
  \quad\text{ for }x\in \mathopen[0, \infty\mathclose[.
\]
If $\alpha,\beta\in\R$ are arbitrary, then~\eqref{Eq:Quot_Phi_incr_bounds}
holds iff  $f_\alpha$ is strictly increasing and $f_\beta$ is strictly decreasing.
Now for $x\in \mathopen]0, \infty\mathclose[$, the derivative $f_\gamma'(x)$
has the same sign as $\gamma - g(x)$ where 
\[
 g(x) &\coloneqq & \frac{ \phi\left(x-\frac{h}2\right) - \phi\left(x+\frac{h}2\right) }
                { x \left(\Phi(x+\frac{h}2) - \Phi(x-\frac{h}2) \right) } %\\ %,
  \,\ = \,\ \mathrm{e}^{-\frac{h^2}8} 
     \cdot \frac{\sinh\left(\frac{xh}2\right)}{\frac{xh}2}
     \cdot       \frac{h\phi(x)}{\Phi(x+\frac{h}2) - \Phi(x-\frac{h}2) }\\
 &=& \mathrm{e}^{-\frac{h^2}8} 
    \frac{\frac2h\int_0^{\frac{h}2}\cosh(xt)\dd t   }
         {\frac2h\int_0^{\frac{h}2}\cosh(xt)\exp(-\frac{t^2}2)  \dd t},
\]
and the unattained supremum and infimum of $g(x)$ over $x\in \mathopen]0, \infty\mathclose[$
are  $\alpha$ and $\beta$ as defined in~\eqref{Eq:Opt_alpha_beta},
by $\exp(-\frac{t^2}2)>\exp(-\frac{h^2}8)$ and by considering $x\rightarrow\infty$,
and by ``Chebyshev's other inequality'' \cite[Chapter IX]{MitPecFink} for the integral of a product of two monotone functions 
applied to yield $\frac2h\int_0^{\frac{h}2}\cosh(xt)\exp(-\frac{t^2}2)  \dd t 
< \frac2h\int_0^{\frac{h}2}\cosh(xt)\dd t\cdot \frac2h \int_0^{\frac{h}2}\exp(-\frac{t^2}2)  \dd t$ 
and by considering $x=0$. This proves our claim except for the inequalities  
in~\eqref{Eq:Opt_alpha_beta}, of which the first one follows from  
\eqref{Eq:Phi_near_zero} and the second  one is trivial  if  $u \coloneqq \frac{h^2}{12}\ge 1$ 
and follows from $\log(1-u)< -u-\frac{u^2}2 <  -u-\frac{u^2}{10}$ otherwise.
\end{proof}

%%%%%%%%%%%%%%%%%%%%%%%%%%%%%%%%%%%%%%%%%%%%%%%%%%%%%%%%%%%%%%%%%%%%%%%%%%%%
%%%%%%%%%%%%%%%%%%%%%%%%%%%%%%%%%%%%%%%%%%%%%%%%%%%%%%%%%%%%%%%%%%%%%%%%%%%%
%%%%%%%%%%%%%%%%%%%%%%%%%%%%%%%%%%%%%%%%%%%%%%%%%%%%%%%%%%%%%%%%%%%%%%%%%%%%
\section{Lemmas on symmetric hypergeometric laws, proof of the main result} \label{Sec:4}
%%%%%%%%%%%%%%%%%%%%%%%%%%%%%%%%%%%%%%%%%%%%%%%%%%%%%%%%%%%%%%%%%%%%%%%%%%%%
%%%%%%%%%%%%%%%%%%%%%%%%%%%%%%%%%%%%%%%%%%%%%%%%%%%%%%%%%%%%%%%%%%%%%%%%%%%%
%%%%%%%%%%%%%%%%%%%%%%%%%%%%%%%%%%%%%%%%%%%%%%%%%%%%%%%%%%%%%%%%%%%%%%%%
To avoid pedantic repetitions of assumptions below, let us agree that
in this section $F,f,n,\sigma,N,\sigma_0, \tau , G$ are in principle fixed 
and as postulated in Theorem~\ref{Thm:Main}, but that we may nevertheless use reduction
arguments as in the proof of Lemma~\ref{Lem:4.2}, where the case of $N=\infty$ is reduced
to the case of $N<\infty$. We have or put 
\[     
  G(s) \,\ =\,\ \Phi\left(\frac{s-\frac{n}{2}}{\tau}\right)
 &\text{ and } & g(s) \,\ \coloneqq\,\ G(s)-G(s-1)\qquad\text{ for }s\in\R, 
\]
and we note the following corollary to Lemma~\ref{Lem:Quot_Phi_incr_bounds}:

%%%%%%%%%%%%%%%%%%%%%%%%%%%%%%%%%%%%%%%%%%%%%%%%%%%%%%%%%%%%%%%%%%%%%%%%%%%%%%%  Lemma 4.1
\begin{Lem}                                                            \label{Lem:g(s+1)/g(s)}
%%%%%%%%%%%%%%%%%%%%%%%%%%%%%%%%%%%%%%%%%%%%%%%%%%%%%%%%%%%%%%%%%%%%%%%%%%%%%%%%%%%%%%%%%%
Let $s\in\mathopen]\frac{n}2,\infty\mathclose[$. Then
\la                                                    \label{Eq:g(s+1)/g(s)}
 \exp\left(-\frac{s-\frac{n}2}{\tau^2}  \right)
  & < &  \frac{g(s+1)}{g(s)} 
 \,\ <\,\ \exp\left( -\left(1-\frac{1}{12\tau^2}\right)\frac{s-\frac{n}2}{\tau^2}   \right).
\al
\end{Lem}       
\begin{proof}Lemma~\ref{Lem:Quot_Phi_incr_bounds} applied to $h\coloneqq\frac1{\tau}$, 
$x\coloneqq (s-\frac{n}2-\frac12)/\tau$, $y\coloneqq (s-\frac{n}2+\frac12)/\tau$. 
\end{proof}
Let us note that $g(s+1)/g(s)$ in~\eqref{Eq:g(s+1)/g(s)} may alternatively be bounded from 
above by $\exp\left(-(s-\frac{n}2- \frac12)/\tau^2 \right)$, as in~\cite[Proof of Lemma~3.1]{HM},
which however appears to be insufficient  for proving Lemma~\ref{Lem:4.3}  below . 

Let $r\coloneqq\frac{N}2$. If $N<\infty$,
then $N$ is even by~\ref{Lem:Hyp_basics}(d), 
$f=\mathrm{h}_{n,r,r}$ by~\ref{Lem:Hyp_basics}(b) and (d), 
hence
\la                                                             \label{Eq:sigma,sigma_0}
 \sigma^2 \,\ =\,\ \frac{n\left(N-n\right)}{4\left(N-1\right)}
 &\text{ and }&  \sigma_0^2 \,\ =\,\ \frac{n\left(N-n\right)}{4 N}
\al
by~\eqref{Eq:Hyp_mean_var_kappa} and~\eqref{Eq:Def_usual_appr_var}, so that in 
particular $\sigma^2>0$ yields $n \in \{1,\ldots, N-1\}$ and thus $N\ge 2$ and  $r\ge 1$,
and we further have 
\begin{equation}
 \sigma_0^2 \,\ \geq\,\ \frac{3}{16} \   \text{ if }  \ n \in \{1,...,N-1\},  
 \qquad 
\sigma_0^2 \,\ \geq\,\ \frac{1}{4} 
\  \text{ if } \  n \in \{2,...,N-2\} \label{SigmaNull},
\end{equation}
by considering $n$ extremal and $N$ minimal.
If $N=\infty$, then $f=\mathrm{b}_{n,1/2}$ and $n\ge 1$.

%%%%%%%%%%%%%%%%%%%%%%%%%%%%%%%%%%%%%%%%%%%%%%%%%%%%%%%%%%%%%%%%%%%%%%%%%%%%%%%%%%%% Lemma 4.2
\begin{Lem}                                                \label{Lem:Bounds_h_etc_n_even_Neu}       
%%%%%%%%%%%%%%%%%%%%%%%%%%%%%%%%%%%%%%%%%%%%%%%%%%%%%%%%%%%%%%%%%%%%%%%%%%%%%%%%%%%%%%%%%%%%%%
Assume $N<\infty$ and $n$ even. Then
\la
  \frac{\sqrt{2}}{4} &<&\sqrt{\frac{(N-2)N^2}{8(N-1)^3}} 
  \,\ \le\,\ \sigma\,f\!\left(\tfrac{n}2\right)
\al
with equality in the second inequality iff $n=2$ or $n=N-2$
\end{Lem}
\begin{proof} We have $n\ge 2$  and hence $N\ge 4$.
Let $a_k\coloneqq\sigma^2\left(\mathrm{H}_{2k,r,r}\right)
\cdot\left(\mathrm{h}_{2k,r,r}(k)\right)^2$
for $k\in \{1,\ldots,r-1\}$. Then, for $k\le r-2$, we have 
\[
 \frac{a_{k+1}}{a_k} &=&  \frac{(k+1)(r-k-1)}{k (r-k)}
 \left( \frac{(r-k)^2(k+\frac12)(k+1)}{(k+1)^2(r-k)(r-k+\frac12)}  \right)^2
  %\text{siehe Dipl.Arb., p.~27}\ldots \\ &=& 
 \,\ = \,\ \frac{1+\frac{1}{4k(k+1)}}{1+\frac{1}{4(r-k)(r-k+1)}}
\]  
and hence $a_{k}\le a_{k+1}$ iff $k(k+1)\le (r-k)(r-k+1)$ iff $k\le \frac{r-1}2$. 
Hence the sequence $(a_k)$ can attain its minimal value only at $k=1$
or at $k=r-1$, and we have in fact 
$\sigma^2\left(\mathrm{H}_{2,r,r}\right)
= \sigma^2\left(\mathrm{H}_{N-2,r,r}\right)= \frac{N-2}{2(N-1)}$
and $\mathrm{h}_{2,r,r}(1) =\mathrm{h}_{N-2,r,r}(r-1)= \frac{N}{2(N-1)}$
and thus $a_1=a_{r-1}=\frac{(N-2)N^2}{8(N-1)^3}$, 
and the latter expression is strictly decreasing in $N\in[4,\infty[$,
as $\left(\log\frac{(x-2)x^2}{(x-1)^3} \right)'= ((x-2)(x-1)x)^{-1}(4-x)<0$
for $x\in\mathopen]4,\infty\mathclose[$ and hence $a_1 \geq \frac{4}{27} \geq \frac{2}{16}$.
\end{proof} 

%%%%%%%%%%%%%%%%%%%%%%%%%%%%%%%%%%%%%%%%%%%%%%%%%%%%%%%%%%%%%%%%%%%%%%%% Lemma 4.3
\begin{Lem}                                        \label{Lem:Bounds_h_etc_n_even}       
%%%%%%%%%%%%%%%%%%%%%%%%%%%%%%%%%%%%%%%%%%%%%%%%%%%%%%%%%%%%%%%%%%%%%%%%%%%%%%%%%%
If $n$ is even and $N<\infty$, then $n=2k$ with $k\in\{1,\ldots,r-1\}$,
$r\ge 2$,
\la                               \label{Eq:sigma_0^2_n_even}   % (22)
 \frac 14 &\le& \sigma_0^2 \,\,\,=\,\,\, \frac{k\,(r-k)}{2r}
   \,\,\, \le\,\,\, \frac{r}8,
\al 
$\sigma =\sqrt{\frac{2r}{2r-1}}\sigma_0$, and 
\la                                          \label{Eq:f(k)_sigma_n_even_lower}    % (23)
 f(k) &>&  \frac1{\sigma_0^{}\sqrt{2\pi}}\exp\left(\frac{23}{192r} - \frac1{16\sigma_0^2}\right)
   \,\ >\,\  \frac1{\sigma_0^{}\sqrt{2\pi}}\exp\left(-\frac1{16\sigma_0^2}\right), \\
 f(k)                                        \label{Eq:f(k)_sigma_n_even_upper}    % (24)
     &<& \frac1{\sigma_0^{}\sqrt{2\pi}}\exp\left(\frac1{8r}-\frac{23}{384\sigma_0^2}\right) \,\ < \,\ \frac{1}{\sigma_0\sqrt{2\pi}}\mathrm{e}^{-\frac{1}{24\sigma_0^2}}
  \,\ <\,\ \frac1{\sigma\sqrt{2\pi}}.
\al
\end{Lem}
\begin{proof} Only the claims in~\eqref{Eq:f(k)_sigma_n_even_lower}
and~\eqref{Eq:f(k)_sigma_n_even_upper} are not obvious.
Writing the binomial coefficient occurring in $f(k)$ in terms of gamma functions
and using the definition of the function $w$ from Lemma~\ref{Lem:w(x)} shows 
that 
\[
 h(k) &\coloneqq &  \log\left(\sigma_0\sqrt{2\pi}f(k)\right)
\]
admits the representation
\[
 h(k) &=& -\log\left(\sqrt{\pi r}w(r)\right) + \log\left(\sqrt{\pi k}w(k)\right) 
  +  \log\left(\sqrt{\pi(r-k)}w(r-k)\right),
\]
so that Lemma~\ref{Lem:w(x)} yields 
\[
 h(k) &>& \frac{23}{192r}-\frac1{8k}-\frac1{8(r-k)}
  \,\ =\,\  \frac{23}{192r} - \frac1{16\sigma_0^2}
  \,\ >\,\ - \frac1{16\sigma_0^2}
\]
and hence~\eqref{Eq:f(k)_sigma_n_even_lower},  
and, using also~\eqref{Eq:sigma_0^2_n_even}, 
\[
  h(k) &<&\frac1{8r}-\frac{23}{192}\left(\frac1k+\frac1{r-k} \right)
  \,\,\,=\,\,\, \frac1{8r}-\frac{23}{384\sigma_0^2} 
  \,\ \le \,\ \left(\frac1{64}-\frac{23}{384}\right)\frac1{\sigma_0^2} \\
  &=& -\frac{17}{384\sigma_0^2} \,\ < \,\ -\frac{1}{24\sigma_0^2} \,\,\, \le\,\,\, - \frac{1}{3r}
\]
and since 
$\log\frac{2r}{2r-1}<\frac1{2r-1} \le \frac2{3r}$ due to $r\ge2$,
we also get
\[
 \frac1{8r} - \frac{23}{384\sigma_0^2} +\log\frac{\sigma}{\sigma_0}
 &<&  -\frac{1}{3r} + \frac1{3r} 
\,\,\,=\,\,\, 0
\]
and hence~\eqref{Eq:f(k)_sigma_n_even_upper}  .
\end{proof}

%%%%%%%%%%%%%%%%%%%%%%%%%%%%%%%%%%%%%%%%%%%%%%%%%%%%%%%%%%%%%%%%%%%%%%%%%%%%%%%% Lemma 4.4
\begin{Lem}                                  \label{Lem:4.2}                     
%%%%%%%%%%%%%%%%%%%%%%%%%%%%%%%%%%%%%%%%%%%%%%%%%%%%%%%%%%%%%%%%%%%%%%%%%%%%%%%%%%
{\rm\textbf{(a)}} 
$f/g$ is strictly decreasing on $\{s\in\Z : \frac{n}2 \le s\le (n\wedge r) +1\}$.  

\smallskip\noindent{\rm\textbf{(b)}} 
We have $f(s)<g(s)$ for $s\in\Z$ with $s>\lfloor n/2 \rfloor$.

\smallskip\noindent{\rm\textbf{(c)}} 
We have $0<F(s)-G(s)<F(\lfloor n/2 \rfloor)-G(\lfloor n/2 \rfloor)$ 
for $s\in\Z$ with $s>\lfloor n/2 \rfloor$.
\end{Lem}
\begin{proof}(a) Let $s\in\Z$ with $\frac{n}2\le s\le n\wedge r$.
Then we have 
\la                   \label{Eq:h(s+1)/h(s)_le_Theta...}
 \quad\frac{f(s+1)}{f(s)}&\le&\Theta\exp\left(-\frac{s-\frac n2}{\sigma^2_0}\right)
 \textstyle
 \quad\text{ with } \Theta\,\coloneqq \,\frac{1-y}{1+y}\mathrm{e}^{2y}\,<\,1
 \text{ where }y\,\coloneqq \,\frac{2s-n+1}{n+1}, 
\al
where $\Theta<1$ holds by Lemma~\ref{Lem:(1+x)/(1+y)} with $x\coloneqq -y$ 
using $y>0$, and the other inequality claimed  holds first  
in case of $r<\infty$, as then Lemma~\ref{Lem:(1+x)/(1+y)} applied
to $\tilde{y}\coloneqq \frac{2s-n+1}{2r-n+1}>0$ and $\tilde{x}\coloneqq -\tilde{y}$ yields,
using $s\ge n-r$ due to $n-r<\frac{n}2$ in the first step,
\[
 \frac{f(s+1)}{f(s)}&=& \frac{(r-s)(n-s)}{(s+1)(r-n+s+1)}
 \,\,\,=\,\,\,\frac{1-y}{1+y} \cdot \frac{1-\tilde{y}}{1+\tilde{y}}
 \,\,\,<\,\,\,\Theta\mathrm{e}^{-2(y+\tilde{y})}
\]
with $2(y+\tilde{y}) = \frac{8(r+1)}{(n+1)(2r-n+1)}(s-\frac{n-1}2)$,  and 
as we have
\[
 \frac{s-\frac n2}{\sigma^2_0}\big/\big(2(y+\tilde{y})\big)
 &=& \frac{r(n+1)(2r-n+1)}{(r+1)n(2r-n)} \cdot\frac{2s-n}{2s-n+1}
 \,\,\,\le\,\,\,1
\]
by using  $(2s-n)/(2s-n+1)\le (2(n\wedge r)-n)/(2(n\wedge r)-n+1)$
and considering separately the cases $n\wedge r = r$
and $n\wedge r = n$, and then also for $r=\infty$,
by taking the limit for $r\rightarrow \infty$ 
in~\eqref{Eq:h(s+1)/h(s)_le_Theta...}.
Now~\eqref{Eq:h(s+1)/h(s)_le_Theta...}  yields,  
using $\tau\ge\sigma_0$ and then
Lemma~\ref{Lem:g(s+1)/g(s)} in case of $s>\frac{n}2$, 
\[
 \frac{f(s+1)}{f(s)}&<& \exp\left(-\frac{s-\frac n2}{\sigma^2_0} \right)
 \,\,\,\le\,\,\,\exp\left(-\frac{s-\frac n2}{\tau^2} \right)
 \,\,\,\le\,\,\, \frac{g(s+1)}{g(s)}
\]
and hence $f(s+1)/g(s+1)<f(s)/g(s)$.

(b) By part~(a) and since $h(s)=0$ for $s>n\wedge r$, we can 
and do assume that $s=\lfloor n/2 \rfloor + 1$. Let us first assume 
that $r<\infty$.

If $n$ is even, then $n=2k$ with $k\in\{1,\ldots,r-1\}$ and $s =k+1$,
and we get 
\[ \textstyle
 g(s)&=& g(k+1)- g(k) 
  \,\,\,=\,\,\,\Phi\left({\textstyle\frac1\tau}\right)-\textstyle\frac12
 \,\,\,\ge\,\,\,\Phi\left({\textstyle\frac1\sigma}\right)-\textstyle\frac12
   \,\,\,>\,\,\, \frac1{\sigma\sqrt{2\pi}}\exp\left(-\frac1{6\sigma^2} \right)
\]
by $\tau\le\sigma$ and Lemma~\ref{Lem:Phi_near_zero}, and, using below  several parts
of Lemma~\ref{Lem:Bounds_h_etc_n_even},  we have
\[
 \frac{f(k+1)}{f(k)}&=& \frac{k\,(r-k)}{(k+1)(r-k+1)} 
 \,\,\,=\,\,\, \frac{k\,(r-k)}{1+r+k\,(r-k) }
 \,\,\,\le\,\,\,  \frac{k\,(r-k)}{\frac32r+k\,(r-k) } 
 \,\,\,=\,\,\, \frac 1{1 +\frac{3}{4\sigma_0^2} }
\]
since $r\ge 2$, and hence 
\[
 \log\frac{f(k+1)}{f(k)}&\le& -\log\left(1 +\frac{3}{4\sigma_0^2}\right)
 \,\,\,<\,\,\,  \frac{-\frac{3}{4\sigma_0^2}}{1 +\frac{3}{4\sigma_0^2}}
 \,\,\,\le\,\,\, \frac{-3}{16\sigma_0^2}
\]
by using $\sigma_0^2\ge \frac 14$ in the last step, so that 
\la  \label{log(h(s)/g(s)_n_even}
 \log\frac{f(s)}{g(s)}
 &=&\log\frac{f(k+1)}{f(k)}  - \log g(s) + \log f(k) \\
 &<& \frac{-3}{16\sigma_0^2} + \frac1{6\sigma^2}   \nonumber
 + \log\left(\sigma\sqrt{2\pi}f(k) \right) 
 \,\ < \,\ -\frac{1}{48\sigma^2} 
\al
using in the final step $\sigma_0\le \sigma$ for the first two terms,
and~\eqref{Eq:f(k)_sigma_n_even_upper}  for the last one.

Let now $n$ be odd. Then $n=2k-1$ with $k\in\{1,\ldots,r\}$ and $s=k$,
and we get 
\la                                \label{Eq:gtau_>_n_odd}
  g(s)&=& \textstyle
   2\left(\Phi\left(\frac1{2\tau}\right) -\frac12 \right) 
   \,\,\,\geq \,\,\, \textstyle
2\left(\Phi\left(\frac1{2\sigma}\right) -\frac12 \right)
  \,\,\,>\,\,\, \frac1{\sigma\sqrt{2\pi}}\exp\left(-\frac1{24\sigma^2} \right)
\al
using $\tau \leq \sigma$ and Lemma~\ref{Lem:Phi_near_zero}.
If $k\in\{1,r\}$, then in either case $\sigma^2=\frac14$ and $h(s)=\frac12$,
and~\eqref{Eq:gtau_>_n_odd} yields
\la                             \label{log(h(s)/g(s)_n=1}  
 g(s)  &\ge& 2 \Phi(1)-1 \,\ =\,\ 0.6827... \,\ >\,\ \frac{1}{2} \,\ =\,\ f(s),
\al
and we now assume that $2\le k\le r-1$. We have
\la                                              \label{Eq:h2k-1.h2k}
 f(s)&=& f(k) \,\ = \,\ \frac{r-k+\frac12}{r-k+1} \mathrm{h}_{2k,r,r} (k) 
\al
and, using $(k-\frac12)/k\ge3/4$ due to $k\ge 2$ for the lower bound
and writing $\sigma_{0,2k}:=\sigma_0(\mathrm{H}_{2k,r,r})$, we get
\la                                      \label{Eq:Quot_var_in}
 \frac{\sigma^2}{\sigma^2_{0,2k}}
   &=&  
 \frac{r(k-\frac12)(r-k+\frac12)}{(r-\frac12)k(r-k)}
  \,\ \in \,\ 
  \left] \frac{3}{4} \, , \, \frac{r-k+\frac12}{r-k} \right]      ,
\al 
and then
\la                   \label{log(h(s)/g(s)_n_odd}
\log\frac{f(s)}{g(s)} 
 &=&\log\left(\sigma^{}_{0,2k}\sqrt{2\pi}\mathrm{h}_{2k,r,r}(k)\right)
   + \log \tfrac{r-k+\frac12}{r-k+1} -\log\left(\sigma\sqrt{2\pi}g(s) \right)
   + \log \tfrac{\sigma}{\sigma^{}_{0,2k}}  \\ \nonumber
 &<&\frac{1}{8r}-\frac{23}{384\sigma_{0,2k}^2}+\log \frac{r-k+\frac12}{r-k+1} 
   +\frac{1}{24\sigma^2}+\frac12\log\frac{r-k+\frac12}{r-k+1} \\ \nonumber
 &<&\left(\frac{1}{24}-\frac{23}{384}\cdot\frac34\right)\frac1{\sigma^2}    
   + \frac{1}{8r} + \frac{1}{4} \log \frac{r-\frac12}{r} 
   + h(r,k)\\ \nonumber
 &<& -\frac{5}{1536\sigma^2}
\al
by using at the second step~\eqref{Eq:f(k)_sigma_n_even_upper}  with $n=2k$,
\eqref{Eq:gtau_>_n_odd}, and \eqref{Eq:Quot_var_in},
at the third step~\eqref{Eq:Quot_var_in}, $k\ge 0$,
and the definition of $h(r,k)$ given below,
and at the final step three applications of $\log(1+x)<x$,
one for $\log((r-\frac12)/r)<-1/(2r)$, and the other two
contained in 
\[
 h(r,k) &\coloneqq & \frac{3}{4 }\log \frac{r-k+\frac12}{r-k+1}
         + \frac{1}{2}\log \frac{r-k+\frac12}{r-k} \\
 &<&  -\frac{3}{4} \cdot \frac{1}{2(r-k+1)} + \frac{1}{2} \cdot \frac{1}{2(r-k)}
 \,\ =\,\ -\frac{r-k-2}{8(r-k-1)(r-k)}
\]
which yields $h(r,k)< 0$ always, namely by the above if $k\le r-2$, and 
by $h(r,r-1)=\tfrac{3}{4} \log\tfrac{3}{4} + \tfrac{1}{2} \log \tfrac{3}{2} 
= \tfrac{1}{4} \log \tfrac{3^33^2}{4^3 2^2} =  \tfrac{1}{4} \log \tfrac{243}{256} <0 $ if $k=r-1$.

By  \eqref{log(h(s)/g(s)_n_even},\eqref{log(h(s)/g(s)_n=1},% 
\eqref{log(h(s)/g(s)_n_odd}, 
there is a constant $c >0$ not depending on  $r,n\in\N$
with $n<2r$ satisfying 
$\log(f(s)/g(s) )\le -c\sigma^{-2}$, and this remains true 
also for  the limit case of $r=\infty$. 

(c) By part (b), $F-G$ is strictly decreasing on 
$\{s\in\Z : s\ge \lfloor n/2\rfloor\}$. Hence we get the second 
inequality claimed and, since $F(s)-G(s)=1-G(s)>0$ for $s\ge n$,
also the first one.
\end{proof}

%%%%%%%%%%%%%%%%%%%%%%%%%%%%%%%%%%%%%%%%%%%%%%%%%%%%%%%%%%%%%%%%%%%%%%%%%%%%  Lemma 4.5
\begin{Lem}                                                             \label{Lem:4.3}      
%%%%%%%%%%%%%%%%%%%%%%%%%%%%%%%%%%%%%%%%%%%%%%%%%%%%%%%%%%%%%%%%%%%%%%%%%%%%%%%%%%%%%%               
Let $\tau\in \mathopen[\sigma_0,\sigma\mathclose]$ and 
$M \coloneqq  \frac{n}2+1+\frac{3}{2}\sigma$.

\smallskip\noindent{\rm\textbf{(a)}}
$f(\cdot-1)/g$ is strictly 
increasing on $\{s\in\Z : \frac{n+1}2 \le s \le M\}$.

\smallskip\noindent{\rm\textbf{(b)}}
We have $g(s)<f(s-1)$ for $s\in\Z$ with
$\lceil n/2 \rceil < s \le M$. 

\smallskip\noindent{\rm\textbf{(c)}}
We have $G(s)-F(s-1)<G(\lceil n/2 \rceil)-F(\lceil n/2 \rceil -1)$
for $s\in\Z$ with $\lceil n/2 \rceil < s \le M$. 
\end{Lem}

\begin{proof}
(a) Let $s\in\Z$ with $\frac{n+1}2\le s \le \frac{n}2 +\frac{3}{2}\sigma$.

If $n=1$, or $N$ is finite and $n=N-1$, then $\sigma =\frac12$ and  hence 
$s=\frac{n+1}2$,
and, using  the unimodality of $\varphi$, we indeed get
$\tfrac{f(s-1)}{g(s)} 
 = \tfrac{1/2}{\Phi(\frac{1}{2\tau})- \Phi(-\frac{1}{2\tau})} 
 < \tfrac{1/2}{\Phi(\frac{3}{2\tau})- \Phi(\frac{1}{2\tau})} 
 = \tfrac{f(s)}{g(s+1)}$.
Hence we can assume $1<n<N-1$ and thus 
$N\ge 4$ and $\sigma_0^2 \geq \tfrac{1}{4}$ in what follows, by~\eqref{SigmaNull}.

Let first also $N<\infty$.
We have $n-r <s\le n\wedge r$, for else  we would have one of 
the inequalities $\frac{n+1}2\le n-r$, $n+1\le \frac{n}2+\frac{3}{2}\sigma$,
$r+1\le \frac{n}2+ \frac{3}{2} \sigma$, which are easily checked to 
be false. 
Hence, putting  $x_1\coloneqq \frac{2s-n}{n}$,  
$y_1\coloneqq -\frac{2s-n-2}{n}$,  $x_2\coloneqq \frac{2s-n}{2r-n}$, and $y_2\coloneqq -\frac{2s-n-2}{2r-n}$, 
we have 
\la
 \frac{f(s)}{f(s-1)}&=&\frac{n-s+1}{s}\cdot\frac{r-s+1}{r-n+s}
   \,\ =\,\ \frac{1+y_1}{1+x_1}\cdot\frac{1+y_2}{1+x_2}  \\ \nonumber
 &\ge& \mathrm{e}_{}^{y_1-x_1+y_2-x_2} 
 \,\ = \,\ \exp\left(-\frac{s-\frac{n+1}2}{\sigma_0^2} \right) 
\al
where the inequality is  a trivial equality  if  $s=\frac{n+1}2$, 
and follows otherwise  by two applications of Lemma~\ref{Lem:(1+x)/(1+y)}, 
as $s\ge\frac{n}2+1$ yields $-y_i\ge 0$ and we also have $x_i-\frac23 x_i^2 \ge -y_i$
for $i\in\{1,2\}$, since 
$y_1+x_1-\frac23 x_1^2 = \frac2n-\frac23\left(\frac{2s-n}{n}\right)^2
\ge\frac2n-\frac23\left(\frac{3\sigma}{n}\right)^2  
\ge \frac2n - \frac{6}{n^2}\cdot\frac{n}{4} \ge 0$
and
$y_2+x_2-\frac23 x_2^2 =\frac{2}{2r-n} -\frac23\left(\frac{2s-n}{2r-n} \right)^2 
\ge \frac{2}{2r-n} -\frac23\left(\frac{3\sigma}{2r-n}\right)^2
\ge \frac{2}{2r-n} -\frac{6}{(2r-n)^2}\cdot\frac{2r-n}4 \ge 0$.
On the other hand,  
putting $x\coloneqq (s-\frac{n+1}2)/\tau$ and applying below Lemma~\ref{Lem:g(s+1)/g(s)}
at the first inequality,$\tau \le \sigma$ and 
$\tau x \le \frac32\sigma-\frac12 \le \frac98\sigma^2
 =\frac98\frac{N}{N-1}\sigma_0^2\le \frac 32 \sigma_0^2$ 
and $\tau^2\ge\sigma_0^2\ge\frac14$ at the second,
and $\frac{\tau^2}{N-1}\le \frac{\sigma^2}{N-1}
=\frac{n(N-n)}{4(N-1)^2}\le  \frac{N^2}{16(N-1)^2}\le \frac19$ 
at the third, we get 
\[
  \exp\left(\frac{s-\frac{n+1}2}{\sigma_0^2} \right) \frac{g(s+1)}{g(s)} 
 &<&\exp\left(\frac{\tau x}{\sigma_0^2}\right) 
    \exp\left(-\left(1-\frac1{12\tau^2}\right)
        \left(\frac{x}\tau+\frac1{2\tau^2}\right)\right)\\ 
 & =&\exp\left( \left( \frac{1}{\sigma_0^2}-\frac1{\tau^2}\right)\tau x
               + \frac{\tau x}{12 \tau^4} - \frac1{2\tau^2} + \frac1{24\tau^4} \right)\\
 &\le&\exp\left( \left( \frac{1}{\sigma_0^2}-\frac1{\sigma^2}\right)\tfrac98\sigma^2
                 + \frac{\frac32\sigma_0^2}{12\tau^2\sigma_0^2} - \frac1{2\tau^2} 
                 + \frac1{6\tau^2}          \right)\\
 &=& \exp\left(\left(\frac{9\tau^2}{8(N-1)}
          -\frac{5}{24}\right)\frac1{\tau^2} \right) \\ 
  &\le&  \exp\left(-\frac{1}{12 \tau^2}\right).
\]  
Thus $ (f(s)/g(s+1)) /(f(s-1)/g(s))\ge
\exp\left(\frac{1}{12 \tau^2}\right)   > 1$, also if $N=\infty$,
hence the claim.

\smallskip
(b) By part (a), we can and do assume that $s= \lceil \tfrac{n}{2} \rceil +1$.
Let first also $r<\infty$.

If $n=2k$ is even, then  $s=k+1$, so that 
Lemma~\ref{Lem:Phi_near_zero}, $\tau \geq \sigma_0$, and 
$\sigma_0^2\ge \frac 14$ from Lemma~\ref{Lem:Bounds_h_etc_n_even} 
yield
\[
 g(s)   &=&\Phi\left(\tfrac{1}{\tau}\right)-\tfrac{1}{2} 
  \,\ \le \,\ \Phi\left(\tfrac{1}{\sigma^{}_0}\right)- \tfrac{1}{2} 
  \,\ < \,\ \frac1{\sigma^{}_0\sqrt{2\pi}}\exp\left(-\frac1{6\sigma_0^2}
    +\frac1{90\sigma_0^4} \right)
 \,\ \le \,\ \frac{ \exp\left(-\frac{11}{90\sigma_0^2}\right)}{\sigma_0\sqrt{2\pi}}
\]
and hence an application of \eqref{Eq:f(k)_sigma_n_even_lower} 
and finally $\sigma_0\le\sigma$ yield
\la
 \log \frac{g(s)}{f(s-1)} 
  &<&  -\frac{11}{90\sigma_0^2} +\frac1{16\sigma_0^2} 
  \,\ =\,\  -\frac{43}{720\sigma_0^2} \,\ \le \,\ -\frac{43}{720\sigma^2}.
\al

Let now $n$ be odd. Then n$=2k-1$  with $k\in\{1,...,r\}$ and $s=k+1$. 
If also $k\leq r-1$ and thus $r\geq 2$, then we have,
using Lemma~\ref{Lem:Incr_Phi_vs_phi}   
with $x\coloneqq h\coloneqq \frac1\tau$ at the first inequality, 
$\tau^2\ge\sigma_0^2\ge 3/16$ by \eqref{SigmaNull} at the second, and $\tau^2 \le \sigma^2$  
at the third, 
\[
 g (s) &=& \Phi\left(\tfrac{3}{2\tau}\right) - \Phi\left(\tfrac{1}{2\tau}\right) 
 \,\ <\,\ \frac{1}{\tau\sqrt{2\pi}} \exp\left( -\frac{13}{24\tau^2}
    + \frac{61}{1440\tau^4} \right) \\
 &\leq&  \frac{1}{\sigma_0\sqrt{2\pi}} \exp\left( -\frac{13}{24\tau^2}
          + \frac{61}{270\tau^2} \right) 
 \,\ = \,\ \frac{1}{\sigma_0\sqrt{2\pi}} \exp\left( -\frac{341}{1080\tau^2} \right) \\
 &\leq& \frac{1}{\sigma_0\sqrt{2\pi}} \exp\left( -\frac{341}{1080\sigma^2} \right)
\]
and further, using the second equality in~\eqref{Eq:h2k-1.h2k} 
and then~\eqref{Eq:f(k)_sigma_n_even_lower} for $\mathrm{h}_{2k,r,r}$,
and writing $\sigma_{0,2k}\coloneqq \sigma_0(\mathrm{H}_{2k,r,r})$,
\[
 f(s-1) &=&   f(k) \,\ =  \,\ \mathrm{h}_{2k,r,r}(k) \cdot \frac{r-k+\frac12}{r-k+1} 
  \,\ > \,\ 
   \frac{\exp \left( 
   - \frac{1}{16\sigma_{0,2k}^2} \right)}{\sigma_{0,2k}\sqrt{2\pi}}    \cdot \frac{r-k+\frac12}{r-k+1} 
\]
and together  with $ \frac{\sigma^2}{\sigma_{0,2k}^2} 
  \leq \frac {r-k+\frac12}{r-k} \leq \frac32$ by \eqref{Eq:Quot_var_in} 
we get the first two inequalities  below and recall $r\ge 2$ for the last: 
\[
 \log \frac{g(s)}{f(s-1)} &\leq& \log \frac{\sigma_{0,2k}}{\sigma_0} - \frac{341}{1080\sigma^2} 
  + \frac{1}{16\sigma_{0,2k}^2}- \log \frac{r-k+\frac12}{r-k+1}\\
&<& -\frac{1}{\sigma^2}\left(\frac{341}{1080} - \frac{3}{32} \right)
 + \frac{1}{2}\log \frac{k(r-k)}{(k-\frac{1}{2})(r-k+\frac{1}{2})} +\log\frac{r-k+1}{r-k+\frac12}\\
&=& -\frac{959}{4320\sigma^2} + \tfrac{1}{2}\log(1+\tfrac{1/2}{k-1/2} ) 
   +  \tfrac{1}{2}\log \left(1+\tfrac{1/2}{r-k+\frac12}\right)
  +\tfrac{1}{2} \log\tfrac{(r-k)\cdot(r-k+1)}{(r-k+\frac12)^2}  \\
&\leq& -\frac{959}{4320\sigma^2} + \frac{1}{2} 
 \left(\frac{1/2}{k-\frac12}+\frac{1/2}{r-k+\frac12}\right)\\
&=& -\frac{1}{\sigma^2}
 \left(\frac{959}{4320}- \frac{1}{8}\cdot \frac{r}{r-\frac12}\right) 
 \,\ \leq \,\ -\frac{1}{\sigma^2} \left(\frac{959}{4320}- \frac{1}{6}\right) 
\,\ = \,\ -\frac{239}{4320\sigma^2}.
\]
%with the help of $\log(1+x) \leq x$. 
If, on the other hand, $k=r$, then $n=2r-1$, $s=r+1$, and hence 
\[
 g(s)&=&G(s)-G(s-1)\,\ =\,\ \Phi\left(\tfrac{3}{2\tau}\right)- \Phi\left(\tfrac{1}{2\tau}\right) 
 \,\ <\,\ 1-\Phi(0) \,\ = \,\  \frac{1}{2}\,\ = \,\ f(s-1).
\]

Hence $g(s)<f(s-1)$ in every case, also if $r=\infty$.

(c) By part (b), $G-F(\cdot-1)$ is strictly decreasing on
$\{s\in\Z: \lceil{n/2}\rceil\le s\le M\}$. 
\end{proof}

%%%%%%%%%%%%%%%%%%%%%%%%%%%%%%%%%%%%%%%%%%%%%%%%%%%%%%%%%%%%%%%%%%%%%%%%%%%
\begin{Lem}                             \label{Lem:4.4}        %  Lemma 4.6
Let $s\in\Z$ with 
$s\ge \frac{n}2+1+ \frac{3}{2}\sigma$. Then
$G(s)-F(s-1) < \frac{1}{\sigma} \phi\left( \frac{3}{2}\right)
 = \frac{0.1295\ldots}{\sigma}$.
\end{Lem}
\begin{proof} 
By Lemma~\ref{Lem:4.2}(c), applicable due to $s-1>\lfloor\frac{n}2\rfloor$,
and then since   $\mathopen[0,\infty\mathclose[ \ni x\mapsto \phi(x)$ and
$\mathopen[1,\infty\mathclose[ \ni x \mapsto x\phi(x)$   
are strictly decreasing and since we have $\frac{3\sigma}{2\tau} \geq \frac{3}{2} \geq 1$,
we get
\begin{align}
 G(s)-F(s-1)&\,\ =\,\  G(s-1)-F(s-1)+G(s)-G(s-1) \nonumber \\
  & \,\ <\,\   G(s)-G(s-1) 
  \,\ < \,\ \tfrac1{\tau}\phi\left(\tfrac{s-\frac{n}2-1}{\tau}  \right) 
 \ \le \  \tfrac1{\tau} \phi\left(  \tfrac{3\sigma}{2\tau}\right)
 \,\ \le\,\  \tfrac1{\sigma} \phi\left(\tfrac{3}{2}\right). \nonumber \qedhere
\end{align}
\end{proof}

%%%%%%%%%%%%%%%%%%%%%%%%%%%%%%%%%%%%%%%%%%%%%%%%%%%%%%%%%%%%%%%%%%%%%%%%%%%%%%%%%%%%%%%
\begin{proof}[Proof of Theorem \ref{Thm:Main}]
%%%%%%%%%%%%%%%%%%%%%%%%%%%%%%%%%%%%%%%%%%%%%%%%%%%%%%%%%%%%%%%%%%%%%%%%%%%%%%%%%%%%%%
Let $d$ be  defined by the first equality in~\eqref{Eq:Def_d}. 
Then the second equality in~\eqref{Eq:Def_d} and
the equality in~\eqref{Eq:d_n_odd_even}  hold by  the first part of 
Lemma~\ref{Lem:Symmetry_and_distances}.

Let us now consider the lower bound
\la            \label{Eq:d_ge_c/sigma}
 d &\ge& \frac{\Phi\left(1\right)-\frac12}{2\sigma}
\al
claimed in~\eqref{Eq:d_in_interval}. If $n$ is odd, then, 
using first $\tau\le\sigma$, and then the concavity of $\Phi$ on $[0,\infty[$ 
and $\sigma\ge\frac12$, we get
\[
 d &=&\Phi\left(\tfrac{1}{2\tau}\right)-\tfrac12 
  \,\ \ge \,\  \frac{\Phi(\tfrac{1}{2\sigma})-\tfrac12}{\tfrac{1}{2\sigma}  } \cdot\tfrac{1}{2\sigma} 
  \,\ \ge \,\  \frac{\Phi(1)-\tfrac12}{2\sigma}
  \,\ = \,\ \frac{0.170672\ldots}{\sigma}
\]   
with equality throughout if $\tau=\sigma$ and $n=1$.
If $n$ is even, then  Lemma~\ref{Lem:Bounds_h_etc_n_even_Neu}  yields 
\[
 d &=& \frac12f\left(\tfrac{n}2\right) \,\ > \,\ \frac{\sqrt{2}}{8\sigma} 
   \,\ = \,\ \frac{0.176776 \ldots}{\sigma}.
\]
The above implies~\eqref{Eq:d_ge_c/sigma} and half of the  
optimality claim in part (b) of the theorem.

We now prove \eqref{Eq:Main_ineq}, using the second part of Lemma~\ref{Lem:Symmetry_and_distances}.
We have~\eqref{Eq:Main_ineq_1} by Lemma~\ref{Lem:4.2}(c). 
To prove~\eqref{Eq:Main_ineq_2}, let $s\in\Z$ with 
$s>\left\lceil\textstyle\frac{n}2\right\rceil$ be given. 
If $s\le M$ from Lemma~\ref{Lem:4.3}, then $G(s)-F(s-1)<d$ by part (c) of that Lemma.
If, on the other hand, $s>M$, then $G(s)-F(s-1)<d$
by Lemma~\ref{Lem:4.4} combined with~\eqref{Eq:d_ge_c/sigma}.  
Hence for  part~(a) of the theorem it only remains to prove the claim involving
the upper bound
\la        \label{Eq:d_<_c/sigma}   
  d &<& \frac{1}{\sigma\sqrt{8\pi}}
\al
contained in~\eqref{Eq:d_in_interval}.
If $n$ is even, then~\eqref{Eq:d_<_c/sigma} follows from~\eqref{Eq:f(k)_sigma_n_even_upper}.

If  $n$ is  odd and $N\neq 2$, then $N\ge 4$.
If $N$ is finite, then, using $\tau\ge \sigma_0$ in the first step,  
$\sigma_0^2\le \frac{N}{16}$ by   \eqref{Eq:sigma_0^2_n_even} and the concavity of $\Phi$ on $[0,\infty[$
and $\frac{\sigma}{\sigma_0}=\sqrt{\frac{N}{N-1}}$ in the second,  
\eqref{Eq:Phi_near_zero} in the third,
$x\coloneqq \frac{1}{N}\in\mathopen]0,\frac 14\mathclose]$ and 
$h(x)\coloneqq -\frac12 \log(1-x) - \frac{2}{3}x + \frac{8}{45}x^2$ in the 
fourth, and the convexity of $h$ on $[0,\frac 14]$
and $h(0)=0$ and 
$ h\left(\frac14\right)  = -\frac12\log\left(\frac34\right)-\frac16+ \frac1{90} < -0.0117145 < 0$ 
in the fifth, we get 
\[
 \sqrt{8\pi}\sigma  \left(\Phi\left(\tfrac{1}{2\tau}\right)-\tfrac12\right) 
  &\le& \frac{\sigma}{\sigma_0} \sqrt{2\pi}
 \frac{\Phi\left(\tfrac{1}{2\sigma_0}\right)-\tfrac12}{\frac1{2\sigma_0}} \\
  &\le& \sqrt{\tfrac{N}{N-1}}\sqrt{2\pi} \frac{\Phi\left(\tfrac{2}{\sqrt{N}}\right)-\tfrac12}{\frac2{\sqrt{N}}}   \\
  &<& \sqrt{\tfrac{N}{N-1}}
    \exp\left( -\tfrac{1}{6 }\left(\tfrac{2}{\sqrt{N}}\right)^2  + \tfrac{1}{90}\left(\tfrac{2}{\sqrt{N}}\right)^4  \right)  \\
  &=& \exp \left(h(x)\right)\\
  &< & 1. 
\]  
If $N=\infty$, then $\tau = \sigma $ and hence $\sqrt{8\pi}\sigma\big( \Phi\left( \frac{1}{2\tau} \right) - \frac 12 \big) <1$ obviously by $\Phi'(x)<\frac1{\sqrt{2\pi}}$ for $0\neq x\in\R$.

Finally if $N=2$, then $n=1$, $\sigma=\frac12$, $\sigma_0=\frac1{\sqrt{8}}$, and
$\sqrt{8\pi}\sigma d = \sqrt{2\pi} \left(\Phi\left(\tfrac{1}{2\tau}\right)-\tfrac12\right) =\varrho(\tau)$
is, as a function of $\rho\in[\sigma_0,\sigma]$, strictly decreasing with 
$\varrho(\sigma_0)=  1.05616\ldots > 1$, and $\varrho(\tau)<1$ iff
$\tau > \tau_0$ with $\tau_0= 0.391961%7
\ldots$,
and we have $\tau_0/\sigma = 0.783923\ldots$.

This proves part (a), and the remaining half of the optimality claim in (b) follows
from 
$\lim_{\sigma\rightarrow\infty}\sigma\,\left(\Phi\left(\frac1{2\sigma}\right)-\frac12 \right)
= \frac1{\sqrt{8\pi}}$. 
\end{proof}
\begin{proof}[Proof of Remark~\ref{Rem:Complements_to_main_result}]
(a) is trivial.

(b) The first inequality in~\eqref{Eq:d-tau-bounds} is trivial by $\frac{N}{N-1}\ge2$ 
and the concavity  of $\Phi$ on $[0,\infty[$.  
If $n$ is odd, then ~\eqref{Eq:d_n_odd_even}, 
concavity again, and $\tau\ge\sigma_0\ge\sqrt{\frac{N-1}{4N}}$
yield the second inequality through
$2\tau  d= \left( \Phi\left(\frac1{\tau}\right) - \frac12 \right)/\left( \frac1{2\tau} \right)
 \ge\left(\Phi\left(\sqrt{\frac{N}{N-1}}\right)-\frac12\right)/\left(\sqrt{\frac{N}{N-1}}\right)$.
If $n$ is even,  and first also $4<N<\infty$, then  
Lemma~\ref{Lem:Bounds_h_etc_n_even_Neu}  yields  
$2\tau d \ge 2\sigma_0 d = 2 \sqrt{\frac{N-1}{N}}\sigma d 
\ge \sqrt{\frac{N-1}{N}}\sqrt{\frac{(N-2)N^2}{8(N-1)^3}} =  \sqrt{\frac{1}{8}(1-\frac{1}{N^2-2N+1})} \geq \sqrt{\frac{1}{8}(1-\frac{1}{36-12+1})} = 0.3461... \geq \Phi \left(1\right) -\frac 12 
\ge \sqrt{\frac{N-1}{N}} \left(\Phi\left(\sqrt{\frac{N}{N-1}}\right)-\frac12\right)$.
If $N=\infty$, then again $\sigma = \sigma_0$  and the claim follows from \eqref{Eq:d_in_interval}. It remains $N=4$, but in this case $2\sigma_0 d = \frac 13 >0.3255... = \sqrt{\frac{3}{4}} \big( \Phi \left( \sqrt{\frac{3}{4}}\right) -\frac 12 \big)$.

If $n$ is even, then $N\neq 2$, and then 
the third inequality 
$d<\frac1{\tau\sqrt{8\pi}}$ follows
trivially from Theorem~\ref{Thm:Main} and $\tau\le\sigma$. 
If $n$ is odd, then $d= \Phi\left(\frac1{2\tau} \right)-\frac12 <\frac1{\tau\sqrt{8\pi}}$. 
  
(c) The second inequality is again obvious. 
In the first inequality, we have equality  if $n$ is odd, 
and if $n$ is even, then Lemma~\ref{Lem:Bounds_h_etc_n_even} and Lemma~\ref{Lem:Phi_near_zero} yield 
\begin{equation}
 d \,\ =\,\ \frac 12 f\left( \frac n2 \right) 
 \,\ \leq \,\ \frac{1}{2\sigma_0\sqrt{2\pi}} \cdot \mathrm{e}^{-\frac{1}{24\sigma_0^2}} 
 \,\ = \,\ \frac{1}{\sqrt{2\pi}} \cdot \frac{1}{2\sigma_0} \mathrm{e}^{-\frac{1}{6(2\sigma_0)^2}} 
 \,\ < \,\ \Phi \left( \frac{1}{2\sigma_0} \right) - \frac 12. \qedhere
\end{equation}
\end{proof}

\addtocontents{toc}{\protect\setcounter{tocdepth}{0}}

\section*{Acknowledgement}% s
We thank Bero Roos for pointing out to us the references \cite{Foley_Hill_Spruill},
\cite{Bobkov_Chist}, and \cite{Mohamed_Mirakhmedov}. We further thank 
Todor Dinev and Christoph Tasto for help with the proofreading.

\addtocontents{toc}{\protect\setcounter{tocdepth}{3}}

\end{document}